\documentclass[final]{siamltex}

\usepackage{latexsym,bm}
\usepackage{mathrsfs,amsmath,amssymb,graphicx} 
\usepackage{amscd} 
\usepackage{relsize,amsmath}

\makeatletter
\newcommand{\biggg}[1]{{\hbox{$\left#1\vbox to 20.5pt{}\right.\n@space$}}}
\newcommand{\Biggg}[1]{{\hbox{$\left#1\vbox to 23.5pt{}\right.\n@space$}}}
\newcommand{\bigggg}[1]{{\hbox{$\left#1\vbox to 26.5pt{}\right.\n@space$}}}
\newcommand{\Bigggg}[1]{{\hbox{$\left#1\vbox to 29.5pt{}\right.\n@space$}}}
\newcommand{\biggggg}[1]{{\hbox{$\left#1\vbox to 32.5pt{}\right.\n@space$}}}
\newcommand{\Biggggg}[1]{{\hbox{$\left#1\vbox to 35.5pt{}\right.\n@space$}}}
\newcommand{\bigggggg}[1]{{\hbox{$\left#1\vbox to 38.5pt{}\right.\n@space$}}}
\newcommand{\Bigggggg}[1]{{\hbox{$\left#1\vbox to 41.5pt{}\right.\n@space$}}}
\makeatother

\newcommand{\mb}[1]{\mathbb #1}

\DeclareMathOperator{\esssup}{ess_{\cdot }  \, \sup}%ess.sup・ｽﾌ為ゑｿｽ
\title{Global asymptotics toward 
rarefaction waves\\ for solutions of
the scalar conservation law\\ 
with nonlinear viscosity}

\author{ 
Akitaka Matsumura \thanks{
Professor Emeritus, 
Department of Pure and Applied Mathematics, 
 Graduate School of Information Science and Technology, 
 Osaka University, Suita, Osaka 565-0871, Japan
({\tt akitaka@math.sci.osaka-u.ac.jp}).}
\and 
Natsumi Yoshida \thanks{
OIC Research Organization, 
Ritsumeikan University, Ibaraki, Osaka 567-8570, Japan
({\tt 14v00067@gst.ritsumei.ac.jp})
/Faculty of Culture and Information Science, 
Doshisha University, Kyotanabe, Kyoto 610-0394, Japan
({\tt jt-bnk68@mail.doshisha.ac.jp}).}}

\begin{document}

\maketitle

\begin{abstract}
In this paper, we investigate the asymptotic behavior of solutions 
to the Cauchy problem 
for the scalar viscous conservation law 
where the far field states are prescribed. 
Especially, we deal with the case 
when 
% the flux function is fully convex, 
% and also 
the viscosity 
is of non-Newtonian type, including a pseudo-plastic case. 
When the corresponding Riemann problem 
for the hyperbolic part 
admits a Riemann solution which 
consists of single rarefaction wave, 
under a condition on nonlinearity of the viscosity,
it is proved that 
the solution of the Cauchy problem tends toward the rarefaction wave 
as time goes to infinity, without any smallness conditions.

%
%An example of SIAM \LaTeX\ macros is presented. Various
%aspects of composing manuscripts for SIAM's journal series
%are illustrated with actual examples from accepted
%manuscripts. SIAM's stylistic standards are adhered to
%throughout, and illustrated.
\end{abstract}



\pagestyle{myheadings}
\thispagestyle{plain}
\markboth{A. MATSUMURA and N. YOSHIDA}{STABILITY OF RAREACTION WAVE}

\section{Introduction}
In this paper, 
we consider the asymptotic behavior of solutions to the Cauchy problem 
for a one-dimensional scalar conservation law with nonlinear viscosity 
\begin{eqnarray}
 \left\{\begin{array}{ll}
  \partial_tu +\partial_x \big( \, f(u) - \sigma(\partial_xu )\, \big) =0
%   = \mu \, 
%     \partial_x \left( \, 
%     \left( \, 1 + ( \, \partial_xu \, )^2 \, \right)^{\frac{p-1}{2}} \partial_xu \, 
%     \right)
  \qquad &(t>0, \: x\in \mathbb{R}), \\[5pt]
  u(0, \, x) = u_0(x) \qquad &( x \in \mathbb{R} ),\\[5pt]
  \displaystyle{\lim_{x\to \pm \infty}} u(t, \, x) =u_{\pm}  
  \qquad &\bigl( t \ge 0 \bigr).   
 \end{array}
 \right.\,
\end{eqnarray}
Here, $u=u(t, \, x)$ is the unknown function of $t>0$ and $x\in \mathbb{R}$, 
the so-called conserved quantity, 
% $$
% Q(u,\: \partial_xu):=f(u) - \sigma(\partial_xu )
% $$ is the total flux 
% (that is, 
the functions $f$ and $-\sigma$ stand for the convective flux and 
viscous/diffusive one, respectively, % depending only on $u$, 
% $\mu$ is the viscosity coefficient, 
% $p$ is the flow behavior index satisfying $p>0$, 
$u_0$ is the initial data, 
and $u_{\pm } \in \mathbb{R}$ 
are the prescribed far field states. 
We suppose that $f$ is a smooth function, 
%satisfying $f(0)=f'(0)=0$
and $\sigma$ is a smooth %odd 
function
satisfying 
\begin{equation}
\sigma(0)=0,\quad \sigma'(v)>0\quad (v\in \mathbb{R}),
\end{equation}
and 
for some $p>0$ %($p$ is the so-called flow behavior index) 
% \begin{equation}
% \begin{cases}
% {}
% \sigma(v) \sim | \, v \,|^p\quad (v\to \infty),
% \qquad \sigma'(v) \sim | \, v \,|^{p-1}\quad (v\to \infty), \\[5pt]
% \sigma(v) \sim -| \, v \,|^p\quad (v\to -\infty),
% \qquad \sigma'(v) \sim -| \, v \,|^{p-1}\quad (v\to -\infty), 
% \end{cases}
% \end{equation}
\begin{equation}
|\,\sigma(v)\,| \sim | \, v \,|^p,
\quad |\,\sigma'(v)\,| \sim | \, v \,|^{p-1}\quad (\,| \, v \,| \to \infty \,).
\end{equation}
A typical example of $\sigma$ in the field of viscous fluid, 
where $u$ corresponds to the fluid velocity, is
\begin{equation}
\sigma(\partial_xu) = \mu \, 
    \left( \, 1 + | \, \partial_xu \, |^2 \, \right)^{\frac{p-1}{2}} \partial_xu
    %\qquad (\mu : \mbox{a positive constant}),
\end{equation}
% where $\mu>0$ is the so-called viscous coefficient, 
where $\mu>0$ is a positive constant, 
which describes a nonlinear
relation between the internal stress $\sigma$ and the deformation velocity $\partial_xu$,
and it is noted that
the cases $p>1$, $p=1$ and $p<1$ physically correspond to 
where the fluid is of dilatant type, Newtonian and pseudo-plastic type,
respectively (see \cite{chh1}, \cite{chh2}, 
\cite{chh-ric}, \cite{jah-str-mul}, \cite{liep-rosh}, 
\cite{ma}, \cite{ma-pr-st}, \cite{soc} and so on).
% We are interested in the stability of the rarefaction wave to (1.1), 
% in particular, the case $p<1$. 
We are interested in the global asymptotics 
for the solution of
(1.1), in particular, the pseudo-plastic case $p<1$, since
there seems no results ever on this case.
First, when $u_- = u_+\,(=:\tilde{u})$, we expect 
the solution globally tends toward the constant 
state $\tilde{u}$ as time goes to infinity.
In fact, we can show the following

\medskip

\noindent
{\bf Theorem 1.1.}\quad{\it
Assume the far field states satisfy $u_- = u_+ (=:\tilde{u})$, 
the viscous flux
$\sigma$, {\rm(1.2)}, {\rm(1.3)}, 
and $p>3/7$. 
Further assume the initial data satisfy
$u_0-\tilde{u} \in H^2$.
Then the Cauchy problem {\rm(1.1)} has a 
unique global in time 
solution $u$ 
satisfying 
\begin{eqnarray*}
\left\{\begin{array}{ll}
u-\tilde{u} \in C^0\cap L^{\infty}
( [\, 0, \: \infty \, ) \, ; H^2 ),\\[5pt]
\partial _x u \in L^2\bigl( \, 0,\: \infty \, ; H^2 \bigr),
\end{array} 
\right.\,
\end{eqnarray*}
and the asymptotic behavior 
$$
\lim _{t \to \infty}\sup_{x\in \mathbb{R}} \, 
\left| \, u(t,x) - \tilde{u} \, \right| = 0.
$$
}

Next, we consider the case 
where the convective flux function $f$ is fully convex, 
that is, 
\begin{equation}
f''(u)>0\quad (u \in \mathbb{R}), 
\end{equation}
and $u_-<u_+$.
Then, since the corresponding Riemann problem (cf. \cite{lax}, \cite{smoller}) 
\begin{eqnarray}
 \left\{\begin{array} {ll}
 \partial _t u + \partial _x \bigl( f(u) \bigr)=0 
 %\, \: \; \qquad 
\quad (t>0,\: x\in \mathbb{R}),\\[5pt]
u(0, \, x)=u_0 ^{\rm{R}} (x)
      := \left\{\begin{array} {ll}
         u_-  & \; (x < 0),\\[5pt]
         u_+  & \; (x > 0)
         \end{array}\right.
 \end{array}
  \right.\,
\end{eqnarray}
turns out to admit a single
rarefaction wave solution, 
we expect that 
the solution of the Cauchy problem (1.1) globally tends toward
the rarefaction wave 
as time goes to infinity. 
Here, the rarefaction wave connecting $u_-$ to $u_+$ 
is given by 
\begin{equation}
u^r \left( \, \frac{x}{t}\: ;\:  u_- ,\:  u_+ \right)
= \left\{
\begin{array}{ll}
  u_-  & \; \bigl(\, x \leq f'(u_-)\,t \, \bigr),\\[7pt]
  \displaystyle{ (f')^{-1}\left( \frac{x}{t}\right) } 
  & \; \bigl(\, f'(u_-)\,t \leq x \leq f'(u_+)\,t\,  \bigr),\\[7pt]
   u_+ & \; \bigl(\, x \geq f'(u_+)\,t \, \bigr).
\end{array}
\right. 
\end{equation} 
Then we can show the following

\medskip

\noindent
% {\bf Theorem 1.2} (Main Theorem I\hspace{-.1em}I){\bf .}\quad{\it
% Assume the far field states satisfy $u_- < u_+$, 
% the convective flux $f \in C^4(\mathbb{R})$ satisfy {\rm(1.4)}, 
% the viscous flux $\sigma \in C^2(\mathbb{R})$ satisfy and {\rm(1.7)}, 
% and the flow behavior index satisfy $p>3/7$. 
% Further assume the initial data satisfy
% $u_0-u_0 ^{\rm{R}} \in L^2$, and
% $\partial _xu_0 \in H^1$. 
% Then the Cauchy problem {\rm(1.1)} has a 
% unique global in time 
% solution $u$ 
% satisfying 
% \begin{eqnarray*}
% \left\{\begin{array}{ll}
% u-u_0 ^{\rm{R}} \in C^0\cap L^{\infty}
% \bigl( \, [\, 0, \: \infty \, ) \, ; L^2 \bigr),\\[5pt]
% \partial _x u \in C^0
%                   \cap L^{\infty}\bigl( \, [\, 0,\: \infty \, ) \, ; H^1 \bigr)
%                   \cap L^2_{\rm{loc}}\bigl( \, 0,\: \infty \, ; H^2 \bigr),
% \end{array} 
% \right.\,
% \end{eqnarray*}
% and the asymptotic behavior 
% $$
% \lim _{t \to \infty}\sup_{x\in \mathbb{R}} \, 
% \left| \, u(t,x) - u^r\left(\, \frac{x}{t}\: ;\: u_-, \, u_+ \right) \, \right| = 0.
% $$
% }
{\bf Theorem 1.2.} \quad{\it
Assume the far field states satisfy $u_- < u_+$, 
the convective flux, $f$ {\rm(1.5)}, 
the viscous flux $\sigma$, {\rm(1.2)}, {\rm(1.3)}, 
and $p>3/7$. 
Further assume the initial data satisfy
$u_0-u_0 ^{\rm{R}} \in L^2$, and
$\partial _xu_0 \in H^1$. 
Then the Cauchy problem {\rm(1.1)} has a 
unique global in time 
solution $u$ 
satisfying 
\begin{eqnarray*}
\left\{\begin{array}{ll}
u-u_0 ^{\rm{R}} \in C^0\cap L^{\infty}
\bigl( \, [\, 0, \: \infty \, ) \, ; L^2 \bigr),\\[5pt]
\partial _x u \in C^0
                  \cap L^{\infty}\bigl( \, [\, 0,\: \infty \, ) \, ; H^1 \bigr)
                  \cap L^2_{\rm{loc}}\bigl( \, 0,\: \infty \, ; H^2 \bigr),
\end{array} 
\right.\,
\end{eqnarray*}
and the asymptotic behavior 
$$
\lim _{t \to \infty}\sup_{x\in \mathbb{R}} \, 
\left| \, u(t,x) - u^r\left(\, \frac{x}{t}\: ;\: u_-, \, u_+ \right) \, \right| = 0.
$$
}

\medskip

It should be emphasized again that as far as the global asymptotic 
stability for either constant states or rarefaction waves,
there have been no results for the case $p<1$ 
(pseudo-plastic type viscosity). 
For the case $p=1$ (Newtonian type viscosity), global
nonlinear stability of both rarefaction wave and viscous shock wave 
were first obtained by
Il'in-Ole{\u\i}nik \cite{ilin-oleinik}. 
For the case $p>1$ 
(dilatant type viscosity), 
when the convective flux satisfies (1.5) 
and viscous flux is even the Ostwald-de Waele type 
% (see \cite{de waele}, \cite{ost}, 
% therefore the viscosity term $\partial_x \bigl( \, 
%   \sigma(\partial_xu) \, \bigr)$ itself is $p$-Laplacian-type), 
($p$-Laplacian type, see \cite{de waele}, \cite{ost}), 
that is, 
\begin{equation}
\sigma(v)= \mu \, 
    \left| \, v \, \right|^{p-1} v, 
\end{equation}
Matsumura-Nishihara \cite{matsu-nishi2} 
proved that if the far field states satisfy $u_- = u_+ (=:\tilde{u})$, 
then the solution globally tends toward the constant state $\tilde{u}$, 
and if $u_- < u_+$, then toward the rarefaction wave. 
Yoshida \cite{yoshida3} also obtained 
the precise time-decay estimates of the solution
toward the constant state and the single rarefaction wave. 
For $p \ge 1$,
it is further considered a case where the flux function $f$ is smooth
% and genuinely nonlinear (that is, $f$ is convex function or concave function)
and convex  
on the whole $\mathbb{R}$ except 
a finite interval $I := (a,\, b) \subset \mathbb{R}$, and 
linearly degenerate on $I$, that is, 
\begin{equation}
\left\{
\begin{array}{ll}
  f''(u) >0 & \; \bigl( \, u \in (-\infty ,\: a\, ]\cup [\, b,\: \infty ) \, \bigr),\\[5pt]
  f''(u) =0 & \; \bigl( \, u \in (a,\, b) \, \bigr).
\end{array}\right.
\end{equation} 
Under the conditions $p \ge 1$, $u_{-}<u_{+}$, (1.8), and (1.9), 
it is proved that the unique global in time solution 
to the Cauchy problem (1.1) globally tends toward 
the multiwave pattern of the combination of 
the viscous contact wave and the rarefaction waves 
as time goes to infinity, where the viscous contact wave
is constructed by the linear heat kernel for $p=1$ 
by Matsumura-Yoshida (\cite{matsumura-yoshida}), 
and also by the Barenblatt-Kompaneec-Zel'dovi{\v{c}} solution 
(see also 
\cite{barenblatt}, \cite{carillo-toscani}, \cite{huang-pan-wang}, \cite{kamin}, 
\cite{vaz1}, \cite{vaz2}, \cite{zel-kom} ) 
of the porous medium equation for $p>1$ 
by Yoshida (\cite{yoshida3}). 
Yoshida (\cite{yoshida1}, \cite{yoshida2}, \cite{yoshida4}) also 
obtained the precise time-decay estimates for these stability results. 
On the other hand, %in the case 
%$u_->u_+$, 
under the Rankine-Hugoniot condition 
\begin{equation}
-s\, (\, u_{+} - u_{-}) + f(u_{+}) - f(u_{-}) = 0, 
%\tag{2.5}
\end{equation}
and Ole{\u\i}nik's shock condition 
\begin{eqnarray}
%f(u) - \big( \, s\, (\, u - u_{\pm}\, ) + f(u_{\pm}) \, \big) < 0 
%\quad \big( \, u \in (u_{+}, u_{-}) \, \big), 
- s\, (\, u - u_{\pm}\, ) + f(u) - f(u_{\pm}) 
\left\{\begin{array}{ll}
< 0 
\quad \big( \, u \in (u_{+}, \: u_{-}) \, \big), \\[5pt]
> 0 
\quad \big( \, u \in (u_{-}, \: u_{+}) \, \big),
 \end{array}
 \right.\,
\end{eqnarray}
%which implies 
%\begin{equation}
%\lambda(u_{-}) \geq s \geq \lambda(u_{+})
%\end{equation}
%where $\lambda:=f'$, 
the local asymptotic stability
of viscous shock waves is proved for $p=1$ by 
Matsumura-Nishihara (\cite{matsu-nishi3}), 
%(not only the asymptotics 
%but also the time-decay estimates), 
and very recently for 
any $p>0$, more generally, for the case where smooth $\sigma$ satisfies 
% \begin{equation}
% \begin{cases}
% {}
% \s \in C^{3}(\mb{R}), \quad \sigma(0) = 0, \\[5pt]
% \sigma'(v) > 0 \quad ( v \in \mb{R} ), \\[5pt]
% \displaystyle{\lim_{v\to \pm \infty}} \sigma(v) =\pm \infty, 
% \end{cases}
% \end{equation}
\begin{equation}
\sigma(0) = 0,\qquad
\sigma'(v) > 0 \quad ( v \in \mb{R} ),\qquad
\displaystyle{\lim_{v\to \pm \infty}} \sigma(v) =\pm \infty, 
\end{equation}
by Yoshida (\cite{yoshida5}), though the global asymptotic 
stability is still open.

\medskip

The proofs of Theorem 1.1 and Theorem 1.2 are given 
by a technical energy method, and a 
Sobolev type inequality motivated by an idea in Kanel' (\cite{kanel}). 
Because the proof of Theorem 1.1 
and the proof of Theorem 1.2 with $p \geq 1$ are much 
easier than that of Theorem 1.2
with $0<p<1$, we only show Theorem 1.2 under the assumption $0<p<1$ in 
the present paper. 

This paper is organized as follows. 
In Section 2, we prepare the basic
properties of the rarefaction wave. 
In Section 3, we reformulate the problem 
in terms of the deviation from 
the asymptotic state. Also,
in order to show the global solution in time and
its asymptotic behavior for the reformulated problem, 
we show the strategy how the local existence and
the {\it a priori} estimates are combined.
In the remaining Section 4, Section 5, and Section 6,
we give the proof of the {\it a priori} estimates
step by step
by using a technical energy method.

%%%%%%%%%%%%%%%%
%%%%%%%%%%%%%%%%
%%%%%%%%%%%%%%%%
%%%%%%%%%%%%%%%%
%%%%%%%%%%%%%%%%
%%%%%%%%%%%%%%%%
\smallskip

{\bf Some Notation.}\quad 
We denote by $C$ generic positive constants unless 
they need to be distinguished. 
In particular, use 
% $C(\alpha,\: \beta,\: \cdots )$ 
% or 
$C_{\alpha,\, \beta,\, \cdots }$ 
when we emphasize the dependency on $\alpha,\: \beta,\: \cdots $.
% Use the symbol ``$<\cdot > $'' as
% $
% < x > := \sqrt{1 + x ^2}. 
% $
For function spaces, 
${L}^p = {L}^p(\mathbb{R})$ and ${H}^k = {H}^k(\mathbb{R})$ 
denote the usual Lebesgue space and 
$k$-th order Sobolev space on the whole space $\mathbb{R}$ 
with norms $||\cdot||_{{L}^p}$ and $||\cdot||_{{H}^k}$, 
respectively. 
% For the weight function $w=w(x) \geq 0$, 
% ${L}^2_{w} = {L}^2_{w}(\mathbb{R})$ and ${H}^k_{w} = {H}^k_{w}(\mathbb{R})$ 
% also denote the weighted Lebesgue space and 
% $k$-th order weighted Sobolev space
% satisfying  
%and 
%the measurable functions $f$ 
%satisfying $\sqrt{w}\, f \in L^2$ and $\sqrt{w}\, f \in {H}^k$ 
%with the norms 
%$$
%f \in L^2_{w} \Longleftrightarrow 
%\sqrt{w}\, f \in L^2, \quad 
%||\, f \, ||_{L^2_{w} }
%:= \left( \, \int_{-\infty}^{\infty} 
%w(x) \, | \, f(x) \,|^2 \, \mathrm{d}x \, \right)^{\frac{1}{2}} < \infty  
%$$
%and 
%$$
%f \in {H}^k_{w} \Longleftrightarrow 
%\sqrt{w}\, \sum _{l=0}^k \, \partial_{x}^l f \in L^2, \quad 
%||\, f \, ||_{{H}^k_{w} }
%:= \left( \, 
%\sum _{l=0}^k \, 
%||\, \partial_{x}^l f \, ||_{L^2_{w} }^2
%\, \right)^{\frac{1}{2}} < \infty,   
%$$
%respectively. 
%If $C^{-1} \leq w(x) \leq C$ for some $C>0$, then we note that 
%$L^2_{w}=L^2$ with $||\cdot||_{{L}^2_{w} }\sim ||\cdot||_{L^2}$
%and ${H}^k_{w}={H}^k$ with $||\cdot||_{{H}^k_{w}}\sim ||\cdot||_{{H}^k}$. 

\medskip

\section{Preliminaries} 
In this section, 
we %shall arrange
prepare a couple of lemmas concerning with 
the basic properties of 
the rarefaction wave.
% for accomplishing the proof of the main theorem. 
Since the rarefaction wave $u^r$ is not smooth enough, 
we need some smooth approximated one. 
% We start with the well-known arguments on $u^r$ 
% and the method of constructing its smooth approximation. 
% We first consider the rarefaction wave solution $w^r$ 
We start with the rarefaction wave solution $w^r$ 
to the Riemann problem 
for the non-viscous Burgers equation:
\begin{equation}
\label{riemann-burgers}
  \left\{\begin{array}{l}
  \partial _t w + 
  \displaystyle{ \partial _x \Big( \, \frac{1}{2} \, w^2 \Big) } = 0 
  \quad \quad \qquad \qquad ( t > 0,\: x\in \mathbb{R}),
%   \, \, \; \; \qquad \quad \qquad ( t > 0,\: x\in \mathbb{R}),
\\[7pt]
  w(0, \, x) = w_0 ^{\rm{R}} ( \, x\: ;\: w_- ,\: w_+)
:= \left\{\begin{array}{ll}
w_+ &\quad (x>0),\\[5pt]
w_- &\quad (x<0),
\end{array}
\right.
\end{array}
  \right.
\end{equation}
where $w_\pm \in \mathbb{R}$ are 
the prescribed far field states
satisfying $w_-<w_+$. 
The unique global weak solution 
$w=w^r\left({x}/{t}\: ;\: w_-,\: w_+\right)$ 
of (\ref{riemann-burgers}) is explicitly given by 
\begin{equation}
\label{rarefaction-burgers}
w^r \Big(\, \frac{x}{t}\: ;\: w_-,\: w_+ \, \Big) %:
= 
  \left\{\begin{array}{ll}
  w_{-} & \bigl(\, x \leq w_{-} \, t \, \bigr),\\[5pt]
  \displaystyle{ \frac{x}{t} } & \bigl(\, w_{-}\, t \leq x \leq w_{+}\, t \, \bigr),\\[5pt]
  w_+ & \bigl(\, x\geq w_{+}\, t \, \bigr).
  \end{array}\right.
\end{equation} 
Next, under the condition 
$f''(u)>0\ (u\in \mathbb{R})$ and $u_-<u_+$, 
the rarefaction wave solution 
$u=u^r\left( {x}/{t}\: ;\: u_-,\: u_+ \, \right)$ 
of the Riemann problem (1.6) 
for hyperbolic conservation law 
is exactly given by 
\begin{equation}
u^r\left( \, \frac{x}{t} \: ; \:  u_-,\: u_+ \, \right) 
= (\lambda)^{-1}\Big( w^r\left( \, \frac{x}{t} \: ; \:  \lambda_-,\: \lambda_+ \, \right)\Big)
\end{equation}
which is nothing but (1.7), 
where $\lambda(u):=f'(u)$ and $\lambda_\pm := \lambda(u_\pm) = f'(u_\pm)$. 
We define a smooth approximation of $w^r( {x}/{t}\: ;\: w_-,\: w_+ \,)$ 
by the unique classical solution 
$$
w=w(\, t, \, x\: ;\: w_-,\: w_+ \,)
\in C^{\infty }\big( \, [\, 0,\: \infty )\times \mathbb{R} \,\big)
$$
to the Cauchy problem for the following %hyperbolic 
non-viscous Burgers equation
\begin{eqnarray}
\label{smoothappm}
\left\{\begin{array}{l}
 \partial _t w 
 + \displaystyle{ \partial _x \Big( \, \frac{1}{2} \, w^2 \, \Big) } =0 
 \, \, \; \; \quad \qquad \qquad \qquad \qquad \qquad  \qquad  \: \,
 (\ t>0,\: x\in \mathbb{R}),\\[7pt]
 w(0,x) 
 = w_0(x) 
 := \displaystyle{ \frac{w_-+w_+}{2} 
    + \frac{w_+-w_-}{2}\,K_{q} \, \int_{0}^{x} \frac{\mathrm{d}y}{(1+y^2)^q} }
 \quad (x\in \mathbb{R}), 
\end{array}
\right.
\end{eqnarray}   
where $K_{q}$ is a positive constant such that 
$$
K_{q} \, \int_{-\infty}^{\infty} \frac{\mathrm{d}y}{(1+y^2)^q} =1 
\quad \bigg( \, q>\frac{1}{2} \, \bigg). 
$$
%Then, we define the approximation for 
By applying the method of characteristics, 
we get the following formula 
\begin{eqnarray}
 \left\{\begin{array} {l}
 w(t, \, x)=w_0\bigl( \, x_0(t, \, x) \, \bigr)=
 \displaystyle{ \frac{\lambda_-+\lambda_+}{2} } 
+ \displaystyle{ \frac{\lambda_+-\lambda_-}{2}
\,K_{q} \, \int_{0}^{x_0(t, \, x) } \frac{\mathrm{d}y}{(1+y^2)^q}  } ,\\[7pt]
 x=x_0(t, \, x)+w_0\bigl( \, x_0(t, \, x) \, \bigr)\,t.
 \end{array}
  \right.\,
\end{eqnarray}
By making use of (2.5) similarly as in \cite{matsu-nishi1}, 
we obtain the properties of 
the smooth approximation $w(t,x : w_-,w_+)$ 
in the next lemma.
%%%%%%%%%%%%%%%%%%%%%%%%%%%%%%%%%%%%%%%%%%%%%%%%%%%%%%%%%%%%%%%%%%%%%%%%%%%%%%%%%%%☆
%%%%%%%%%%%%%%%                     w(t,x)                    %%%%%%%%%%%%%%%%%%%%%☆
%%%%%%%%%%%%%%%%%%%%%%%%%%%%%%%%%%%%%%%%%%%%%%%%%%%%%%%%%%%%%%%%%%%%%%%%%%%%%%%%%%%☆

\medskip

\noindent
{\bf Lemma 2.1.}\quad{\it
% Assume that the far field states satisfy $w_-<w_+$.
Assume $w_-<w_+$.  
Then the classical solution $w=w(t, x : w_-,w_+)$
given by {\rm(2.4)} 
satisfies the following properties: 

\noindent
{\rm (1)}\ \ $w_- < w(t, \, x) < w_+$ and\ \ $\partial_xw(t, \, x) > 0$  
\quad  $(t>0, \: x\in \mathbb{R})$.

\smallskip

\noindent
{\rm (2)}\ For any $q>1/2$ and $r \in [\, 1, \: \infty \,]$, there exists a positive 
constant $C_{q,\,r}$ such that
             \begin{eqnarray*}
                 \begin{array}{l}
                    \parallel \partial_x w(t)\parallel_{L^r} \leq 
                    C_{q,\,r} \, (1+t)^{-1+\frac{1}{r}} 
                    \quad \bigl(t\ge 0 \bigr),\\[5pt]
                    \parallel \partial_x^2 w(t) \parallel_{L^r} \leq 
                    C_{q,\,r} \, (1+t)^{-1-\frac{1}{2q} \, 
                    \left( 1 - \frac{1}{r} \right)} 
                    \quad \bigl(t\ge 0 \bigr),\\[5pt]
                    \parallel \partial_x^3 w(t) \parallel_{L^r} \leq 
                    C_{q,\,r} \, (1+t)^{-1-\frac{1}{2q} \, 
                    \left( 2 - \frac{1}{r} \right)} 
                    \quad \bigl(t\ge 0 \bigr).
                    \end{array}       
              \end{eqnarray*}
              
\smallskip

\noindent
{\rm (3)}\; $\displaystyle{\lim_{t\to \infty} 
\sup_{x\in \mathbb{R}}
\left| \,w(t, \, x)- w^r \left( \frac{x}{t} \right) \, \right| = 0}.$
}

\bigskip
%%%%%%%%%%%%%%%%%%%%%%%%%%%%%%%%%%%%%%%%%%%%%%%%%%%%%%%%%%%%%%%%%%%%%%%%%%%%%%%%%%%☆
%%%%%%%%%%%%%%%                     w(t,x)                    %%%%%%%%%%%%%%%%%%%%%☆
%%%%%%%%%%%%%%%%%%%%%%%%%%%%%%%%%%%%%%%%%%%%%%%%%%%%%%%%%%%%%%%%%%%%%%%%%%%%%%%%%%%☆

%\noindent
We now define the approximation for 
the rarefaction wave $u^r\left( {x}/{t}\: ;\: u_-,\: u_+ \, \right)$ by 
\begin{equation}
U(\, t, \, x\: ; \: u_-,\: u_+ \,) 
%:
= (\lambda)^{-1} \bigl( \,w(\, t, \, x\: ;\: \lambda_-,\: \lambda_+\,) \, \bigr).
\end{equation}
Noting the assumption of the smooth flux function $f$, %to be 
%$\lambda'(u)\left( =\frac{\mathrm{d}^2f}{\mathrm{d}u^2}(u)\right)>0$.
we have the next lemma. 

\medskip

\noindent
{\bf Lemma 2.2.}\quad{\it % <----- Lemma 2.2☆☆☆☆☆☆☆☆☆☆☆☆☆☆☆☆
% Assume that the far field states satisfy $u_-<u_+$, 
% and the flux function $f\in C^3(\mathbb{R})$, $f''(u)>0 \; (\,u\in [\,u_-,\:u_+\,]\,)$. 
% Assume that $u_-<u_+$, $f\in C^4(\mathbb{R})$, and $f''(u)>0 \; (\,u\in [\,u_-,\:u_+\,]\,)$.
Assume $u_-<u_+$ and $f''(u)>0 \ (u \in \mathbb{R})$.  
Then we have the following: %properties:

\noindent
{\rm (1)}\ $U(t, \, x)$ defined by {\rm (2.6)} is %<----- 2.6☆☆☆☆☆☆☆☆☆☆☆☆☆☆☆☆
% the unique $C^3$-global solution in space-time 
% of the Cauchy problem
the unique smooth global solution to
the Cauchy problem
%\begin{equation*}
$$
\left\{
\begin{array}{l} 
\partial _t U +\partial _x \bigl( f(U ) \bigr) = 0 \quad
% \, \, \, \, \; \; \; \quad \quad \qquad \qquad \qquad \qquad  
(t>0, \: x\in \mathbb{R}),\\[7pt]
U(0,x) 
= \displaystyle{ (\lambda)^{-1} \left( \, \frac{\lambda_- + \lambda_+}{2} 
+ \frac{\lambda_+ - \lambda_-}{2} 
  \,K_{q} \, \int_{0}^{x} \frac{\mathrm{d}y}{(1+y^2)^q} \, \right) } 
\quad( x\in \mathbb{R}),\\[10pt]
\displaystyle{\lim_{x\to \pm \infty}} U(t, \, x) =u_{\pm} \quad
% \, \, \: \: \; \; \quad \quad \qquad \qquad \qquad \qquad \qquad \qquad 
\bigl(t\ge 0 \bigr).
\end{array}
\right.\,     
$$
%\end{equation*}
{\rm (2)}\ \ $u_- < U(t, \, x) < u_+$ and\ \ $\partial_xU(t,x) > 0$  
\quad  $(t>0, \: x\in \mathbb{R})$.

\smallskip

\noindent
{\rm (3)}\ For any $q>1/2$ and $r \in [\,1,\: \infty \,]$, there exists a positive 
constant $C_{q,\,r}$ such that
             \begin{eqnarray*}
                 \begin{array}{l}
                    \parallel \partial_x %U(t,\cdot \: )
                    U(t) 
                    \parallel_{L^r} \leq 
                    C_{q,\,r}\, (1+t)^{-1+\frac{1}{r}} 
                    \quad \bigl(t\ge 0 \bigr),\\[5pt]
                    \parallel \partial_x^2 U(t) \parallel_{L^r} \leq 
                    C_{q,\,r}\, (1+t)^{-1-\frac{1}{2q} \, 
                    \left( 1 - \frac{1}{r} \right)}
                    \quad \bigl(t\ge 0 \bigr),\\[5pt]
                    \parallel \partial_x^3 U(t) \parallel_{L^r} \leq 
                    C_{q,\,r}\, (1+t)^{-1-\frac{1}{2q} \, 
                    \left( 2 - \frac{1}{r} \right)}
                    \quad \bigl(t\ge 0 \bigr).
                    \end{array}       
              \end{eqnarray*}
              
\smallskip

\noindent
{\rm (4)}\; $\displaystyle{\lim_{t\to \infty} 
\sup_{x\in \mathbb{R}}
\left| \,U(t, \, x)- u^r \left( \frac{x}{t} \right) \, \right| = 0}.$

% \smallskip

% \noindent
% {\rm (5)}\ For any $\epsilon \in (0,\:1)$, there exists a positive 
% constant $C_\epsilon$ such that
% $$
% \left|\, 
% U(t,x)-u_+ \, \right|
% \leq C_\epsilon \, (1+t)^{-1+\epsilon}
%      \mathrm{e}^{-\epsilon \, | \, x-\lambda_+ \, t \, |}
% \quad \bigl(\, t\ge 0, \: x \ge \lambda_+ \, t \, \bigr).
% $$

% \smallskip

% \noindent
% {\rm (6)}\ For any $\epsilon \in (0,\:1)$, there exists a positive 
% constant $C_\epsilon$ such that
% $$
% \left|\, 
% U(t, \, x)-u_- \, \right|
% \leq C_\epsilon \, (1+t)^{-1+\epsilon}
%      \mathrm{e}^{-\epsilon \, | \, x-\lambda_-\, t \, |}
% \quad \bigl(\, t\ge 0, \: x \le \lambda_-\, t\, \bigr).
% $$

% \smallskip

% \noindent
% {\rm (7)}\ For any $\epsilon \in (0,\:1)$, there exists a positive 
% constant $C_\epsilon$ such that
% $$
% \left| \,U(t, \, x) - u^r\left( \frac{x}{t}\right) \, \right| 
% \leq C_\epsilon \, (1+t)^{-1+\epsilon} 
% \qquad \bigl(\, t \ge 1, \: \lambda_-\, t \le x \le \lambda_+\, t\, \bigr).
% $$

% \smallskip

% \noindent
% {\rm (8)}\ For any $\epsilon \in (0,1)$ and $1\le r \le \infty$, 
% there exists a positive 
% constant $C_{\epsilon,\, r}$ such that
% $$ 
% \left|\left|
%  \,U(t,\cdot \: ) - u^r\left( \frac{\cdot }{t}\right) \, 
% \right|\right|_{L^r}  
% \leq C_{\epsilon,\, r}\, (1+t)^{-1+\frac{1}{r}+\epsilon} 
% \qquad \bigl(t \ge 0 \bigr).
% $$
}

\medskip

Because the proofs of them are well-known, 
we omit the proofs here
(see \cite{hashimoto-matsumura}, \cite{hattori-nishihara}, \cite{liu-matsumura-nishihara}, \cite{matsu-nishi1}, 
\cite{matsumura-yoshida}, \cite{yoshida1}, and so on).

\section{Reformulation of the problem} 
In this section, we reformulate our problem (1.1) 
in terms of the deviation from the asymptotic state. 
% At first, note that
% without loss of generality, 
% we can further assume 
% \begin{equation}
% f(0)=f'(0)=0, \quad u_-<0<u_+. 
% \end{equation}
% by changing the variables and constants as
% $x-\tilde{\lambda} \,t \mapsto x$,  $u-b \mapsto u$,
% $f(u+b)-f'(b)\,u-f(a) \mapsto f(u)$ and $a-b \mapsto a$
% in this order. 
% \begin{equation}
% f(0)=f'(0)=0.
% \end{equation}
% In fact, we may
% change the variables 
% $x-f'(0)t \mapsto x$, and  
% $f(u)-f'(0)u-f(0) \mapsto f(u)$. 
Now letting 
\begin{equation}
u(t, \, x) = U(t, \, x) + \phi(t, \, x), 
\end{equation}
we reformulate the problem (1.1) in terms of 
the deviation $\phi $ from $U$ as 
\begin{eqnarray}
 \left\{\begin{array}{ll}
  \partial _t\phi + \partial_x \big( \, f(U+\phi) - f(U) \, \big) \\[5pt]
  \qquad - \partial_x 
    \big( \, 
    \sigma ( \, \partial_x U + \partial_x \phi \, ) 
    - \sigma ( \, \partial_x U  \, )  \, 
    \big)
    = 
    \partial_x \big( \, \sigma ( \, \partial_x U  \,)  \, \big) \quad (t>0,\: x\in \mathbb{R}),\\[5pt]
%     \\[2pt]
%     \quad \quad \qquad \qquad \qquad \qquad \qquad 
%     \qquad \qquad \qquad \; 
%   \, \, \: \: \quad  (t>0,\: x\in \mathbb{R}), \\[5pt]
  \phi(0, \, x) = \phi_0(x) 
  := u_0(x)-U(0, \, x) \quad (x\in \mathbb{R}),
%   \qquad \qquad \quad \; \: \: \; \; \: \: \; \; \,\, \, (x\in \mathbb{R}),
  \\[5pt]
  \displaystyle{\lim_{x \to \pm \infty}} \phi (t,x) =0 \quad  \bigl( t \ge 0 \bigr).
%   \qquad \qquad \qquad \; \: \: \: \: 
%   \qquad \qquad \quad \quad \quad \; \; \; \; \: \: \:\, 
%   \bigl( t \ge 0 \bigr).
 \end{array}
 \right.\,
\end{eqnarray}
Then we look for 
the unique global in time solution 
$\phi $ which has the asymptotic behavior 
\begin{equation}
\displaystyle{
\sup_{x \in \mathbb{R}}\left|\,  \phi(t, \, x) \, \right|
\to  0 \qquad (t\to  \infty)}. 
\end{equation}
Here we note %the fact 
that $\phi_0 \in H^2$ by the assumptions on $u_0$, 
and Lemma 2.2. 
Then the corresponding theorem for $\phi$ to Theorem 1.2  we should prove is 
as follows. 

\medskip

\noindent
{\bf Theorem 3.1.} \quad{\it
% Assume the far field states $u_{\pm}$ and 
% the convective flux $f \in C^4(\mathbb{R})$ satisfy 
% {\rm(1.4)} and {\rm(3.1)}, 
% the viscous flux $\sigma \in C^2(\mathbb{R})$ satisfy {\rm(1.7)}, 
% and the flow behavior index satisfy $p>3/7$. 
% Further assume the initial data satisfy
% $\phi_0 \in H^2$. 
Assume the far field states satisfy $u_- < u_+$, 
the convective flux $f$, {\rm(1.5)}, 
the viscous flux $\sigma$, {\rm(1.2)}, {\rm(1.3)}, 
and $p>3/7$. 
Further assume the initial data satisfy
$\phi_0 \in H^2$. 
Then the Cauchy problem {\rm(3.2)} has a 
unique global in time 
solution $u$ 
satisfying 
\begin{eqnarray*}
\left\{\begin{array}{ll}
\phi \in C^0\cap L^{\infty}
\bigl(\,[\, 0, \: \infty \, ) \, ; H^2 \bigr),\\[5pt]
\partial _x \phi  \in L^2\bigl( \, 0,\: \infty \, ; H^2 \bigr),
\end{array} 
\right.\,
\end{eqnarray*}
and the asymptotic behavior 
$$
\lim _{t \to \infty}\sup_{x\in \mathbb{R}}
|\,\phi (t, \, x)\,| = 0. 
$$
}

% To accomplish the proof of Theorem 3.1, 
% we prepare the local existence precisely, 
% we formulate 
% the problem (3.3) at general 
% initial time $\tau \ge 0$: 
Theorem 3.1 is shown by combining the 
local existence of the solution together with
the {\it a priori} estimates as in the previous papers.
% To state the local existence precisely, 
% we generalize the Cauchy problem
% (3.3) to one at general 
% initial time $\tau \ge 0$: 
To state the local existence precisely, 
the Cauchy problem at general 
initial time $\tau \ge 0$ with the given initial 
data $\phi_\tau \in H^2$ is formulated:
\begin{eqnarray}
 \left\{\begin{array}{ll}
  \partial _t\phi + \partial_x \left( \, f(U+\phi) - f(U) \, \right) \\[5pt]
  \qquad - \partial_x 
    \big( \, 
    \sigma ( \, \partial_x U + \partial_x \phi \,) 
    - \sigma ( \, \partial_x U  \,)  \, 
    \big)
    = 
    \partial_x \big( \, \sigma( \, \partial_x U  \,)  \, \big) \quad
(t>\tau,\: x\in \mathbb{R}),
    \\[5pt]
%     \quad \quad \qquad \qquad \qquad \qquad \qquad 
%     \qquad \qquad \qquad \; 
%   \, \, \: \: \quad  (t>\tau,\: x\in \mathbb{R}), \\[5pt]
  \phi(\tau, \, x) = \phi_\tau(x) 
  := u_\tau(x)-U(\tau, \, x) \quad (x\in \mathbb{R}),
%   \qquad \qquad \quad \; \: \: \; \; \: \: \; \; \,\, \, (x\in \mathbb{R}),
  \\[5pt]
  \displaystyle{\lim_{x \to \pm \infty}} \phi (t, \, x) =0 
  \quad %\qquad \qquad \; \: \: \: \: 
%  \qquad \qquad \quad \quad \quad \; \; \; \; \: \: \:\, 
  \bigl( t \ge \tau \bigr).
 \end{array}
 \right.\,
\end{eqnarray}

\medskip

\noindent
{\bf Theorem 3.2} (local existence){\bf .}\quad{\it
For any $ M > 0 $, there exists a positive constant 
$t_0=t_0(M)$ not depending on $\tau$ 
such that if $\phi_{\tau} \in H^2$ and 
$\displaystyle{ 
\| \, \phi_{\tau} \, \|_{{H}^{2}} \leq M}$, then
the Cauchy problem {\rm (3.4)} has a unique solution $\phi$ 
on the time interval $[\, \tau, \: \tau+t_{0}(M)\, ]$ satisfying 
% \begin{eqnarray*}
% \left\{\begin{array}{ll}
% \phi \in C^0\cap L^{\infty}
% \bigl( \, [\, \tau, \: \tau+t_{0}\, ] \, ; L^2 \bigr),\\[5pt]
% \partial _x \phi  \in C^0
%                   \cap L^{\infty}\bigl( \, [\, \tau, \: \tau+t_{0}\, ] \, ; H^1 \bigr)
%                   \cap L^2\bigl( \, \tau, \: \tau+t_{0} \, ; H^2 \bigr),
% \end{array} 
% \right.\,
% \end{eqnarray*}
\begin{eqnarray*}
\phi \in C^0\big( \, [\, \tau, \: \tau+t_{0}\, ] \, ; H^2 \big) \cap
L^2( \, \tau, \: \tau+t_{0} \, ; H^3 ).
\end{eqnarray*}
}

\medskip

\noindent
The proof of Theorem 3.2 is given by standard iterative method with the aid of
the semigroup theory 
by Kato \cite{kato1}, \cite{kato2}. 
Because the proof is similar to the one in Yoshida \cite{yoshida5}, 
we omit the details here (cf. \cite{lad-sol-ura}, \cite{lions}, \cite{yoshida3}). 
The {\it a priori estimates} we establish in Section 4, Section 5 and Section 6 
are the following. 

\medskip

\noindent
{\bf Theorem 3.3} ({\it a priori} estimates){\bf .}\quad{\it 
Under the same assumptions in Theorem 3.1,
for any initial data 
$\phi_0 \in H^2$, 
% there exists a positive constant $C$
% depending only on $u_{\pm}$, $\| \, \phi_{0} \, \|_{{H}^{2}}$, 
% and the shape of the convective flux $f$ and the viscous flux $\sigma$
there exists a positive constant $C_{\phi_0}$
such that 
if the Cauchy problem {\rm (3.2)}
has a solution 
$\phi$ 
on a time interval $[\, 0, \: T\, ]$ satisfying 
\begin{eqnarray*}
\phi \in C^0\big( \, [\, 0, \: T\, ] \, ; H^2 \big) \cap
L^2( \, 0, \: T \, ; H^3)
\end{eqnarray*}
for some constant $T>0$, 
then it holds that 
\begin{align}
\begin{aligned}
\| \, \phi(t) \, \|_{{H}^{2}}^{2} 
& + \int^{t}_{0} 
   \big\| \, 
   (\, \sqrt{\, \partial_x U} \:  \phi \,)(\tau) 
   \, \big\|_{L^2}^{2} \, \mathrm{d}\tau 
\\
& 
+ \int^{t}_{0}
   \big( 
 \left\| \, \partial_{x}\phi(\tau) \, \right\|_{{H}^2}^{2} 
 + \left\| \, \partial_{t}\partial_{x}\phi(\tau) \, \right\|_{L^2}^{2} 
   \, \big) \, \mathrm{d}\tau 
\leq C_{\phi_{0}}
\qquad \big( \, t \in [\, 0, \: T\, ] \,\big).
\end{aligned}
\end{align}
}

\noindent
% Combining Theorem 3.2 and Theorem 3.3, we have Theorem 3.1.
Once Theorem 3.3 is established,
by combining the local existence Theorem 3.2
with  $M=M_0:=\sqrt{ \, C_{\phi_0}}$,
$\tau = n \, t_0(M_0)$, 
and $\phi_\tau=\phi\big( \, n \, t_0(M_0) \, \big)\ (n=0,1,2,\dots)$ together 
with the {\it a priori} estimates with $T=(n+1) \, t_0(M_0)$
inductively,
the
unique solution of (3.3) 
$\phi 
\in C^0\big( \, [\, 0,\: n\, t_0(M_0)\, ] \, ;H^2\big)
\cap L^2\big( \,0,\: n\, t_0(M_0)\, ;H^3\big)$
for any $n \in \mathbb{N}$
is easily constructed,
that is, the global solution in time
$\phi\in C^0\big( \, [\, 0,\: \infty) \, ;H^2\big)\cap L^2_{\rm{loc}}( \, 0,\: \infty \, ; H^3)$. 
Then, the {\it a priori} estimates again
assert that
\begin{equation}
\sup_{t\ge 0}\, \| \, \phi (t)\,  \|_{H^2}< \infty,\quad
\int _0^\infty \big(\,\| \, \partial_x\phi(t) \, \|_{H^2}^2+
\| \, \partial_t\partial_x\phi(t) \, \|^2_{L^2} \, \big)
\, \mathrm{d}t 
< \infty,
\end{equation}
which easily gives
\begin{align}
\begin{aligned}
\int _0^{\infty }\bigg|\,\frac{\mathrm{d}}{\mathrm{d}t} 
\|  \, \partial_x\phi(t)  \, \|_{L^2}^2 \,\bigg| \, \mathrm{d}t
< \infty.
\end{aligned}
\end{align}
Hence, it follows from (3.6) and (3.7) that
\begin{equation*}
\|  \, \partial_x\phi(t)  \, \|_{L^2} \to 0\quad (t \to \infty).
\end{equation*}
Due to the Sobolev inequality, the
desired asymptotic behavior in Theorem 3.1 is obtained as
\begin{equation*}
\sup_{x \in \mathbb{R}}|\,\phi(t,x)\,|
\leq \sqrt{2}\, \| \,\phi(t)\, \|^{\frac{1}{2}}_{L^2} 
\| \,\partial_x\phi(t)\, \| ^{\frac{1}{2}}_{L^2}
\to 0 \quad (t \to \infty).
\end{equation*} 
Thus, Theorem 3.1 is shown by combining Theorem 3.2 together with
Theorem 3.3. In the following sections, we give the proof of 
the {\it a priori} estimates, Theorem 3.3.
To do that, in the whole remaining sections we assume 
$\phi \in C^0\big( \, [\, 0, \: T\, ] \, ; H^2 \big) \cap
L^2( \, 0, \: T \, ; H^3)$ is a solution of (3.2) for
some $T>0$, and for simplicity we use the notation $C_0$ to denote 
positive
constants which may depend on the initial data $\phi_0 \in H^2$,
and the shape of the equation but not depend on $T$.

\bigskip 

\noindent
%%%%%%%%%%%%%%%%%%%%%%%%%%%%%%%%%%%%%%%%%%%%%%%%%%%%%%%%%%%%%%%%%%%%%%%%%%%%%
\section{A priori estimates I}
In this section, 
we show the following basic $L^2$-energy estimate for $\phi$.

\medskip

\noindent
{\bf Proposition 4.1.}\quad {\it
For $0<p<1$, there exists a positive constant $C_0$
such that 
% \begin{align*}
% \begin{aligned}
% & \| \, \phi(t) \, \|_{{L}^{2}}^{2} 
%  + \int^{t}_{0} 
%    \Big|\Big| \, 
%    \Big(\, \sqrt{\partial_x U} \:  \phi \, \Big)(\tau) 
%    \, \Big|\Big|_{L^2}^{2} \, \mathrm{d}\tau 
%    + \int^{t}_{0} \int ^{\infty }_{-\infty } 
%      <\, \partial_x \phi \, >^{p-1} 
% |\, \partial_x \phi \, |^2 
% \, \mathrm{d}x\mathrm{d}\tau \\
% &+ \int^{t}_{0} \Big( \, 
%    \left\| \, \partial_x\phi(\tau) \, \right\|_{L^{2}(\{ \, |\partial_x \phi|<1 \, \})}^{2} 
%    + \Vert \, \partial_x \phi(\tau) \, 
%        \Vert_{L^{p+1}(\{ \, |\partial_x \phi|\geq1 \, \})}^{p+1}
%    \, \Big) \, \mathrm{d}\tau \\
% & \qquad \qquad \qquad \qquad \qquad \qquad 
% &\leq C \, 
%      \Bigl( \, 
%      \| \, \phi_{0} \, \|_{{L}^2}^2 
%      \, \Bigr) 
% \quad \big( \, t \in [\, 0, \: T\, ] \, \big). 
% \end{aligned}
% \end{align*}
\begin{align*}
\begin{aligned}
\| \, \phi(t) \, \|_{{L}^{2}}^{2} 
 &+ \int^{t}_{0} \int
   |\phi|^2 \,{\partial_x U}\, \mathrm{d}x\mathrm{d}\tau \\
&+ \int^{t}_{0} \int
     <\partial_x \phi>^{p-1} 
|\, \partial_x \phi \, |^2 
\, \mathrm{d}x\mathrm{d}\tau \leq C_0
% \\
% &+ \int^{t}_{0} \Big( \, 
%    \left\| \, \partial_x\phi(\tau) \, \right\|_{L^{2}(\{ \, |\partial_x \phi|<1 \, \})}^{2} 
%    + \Vert \, \partial_x \phi(\tau) \, 
%        \Vert_{L^{p+1}(\{ \, |\partial_x \phi|\geq1 \, \})}^{p+1}
%    \, \Big) \, \mathrm{d}\tau \\
%& \qquad \qquad \qquad \qquad \qquad \qquad 
% &\leq C \, 
%      \Bigl( \, 
%      \| \, \phi_{0} \, \|_{{L}^2}^2 
%      \, \Bigr) 
\qquad \big( \, t \in [\, 0, \: T\, ] \, \big),
\end{aligned}
\end{align*}
where $<s>:= (1+s^2)^{1/2}\ (s \in \mathbb{R})$.
}

\medskip

To obtain Proposition 4.1, 
we first show the uniform boundedness of $\|\phi\|_{L^{\infty}}$ by using 
the $L^q\,(q\ge 2)$ energy estimates as follows 
(cf. \cite{ilin-kalashnikov-oleinik}, \cite{lad-sol-ura}, \cite{lions}). 

\medskip

\noindent
{\bf Lemma 4.1.} 
% (uniform boundedness for $\phi$){\bf .}
\quad {\it
There exists a positive constant $C_0$ such that 
$$
% \displaystyle{ 
% \sup_{x \in \mathbb{R}} \, |\, \phi(t, \,x) \,| \leq C_0
\|\, \phi(t) \,\|_{L^{\infty}} \leq C_0
\qquad \big( \, t \in [\, 0, \: T\, ] \, \big).
$$
}

\medskip

{\bf Proof of Lemma 4.1.}
For $ r \geq 1$, multiplying the equation in (3.2) by 
$|\, \phi \, |^{r-1}\, \phi$, and integrating the resultant formula
with respect to $x$, we have, after integration by parts,
% we obtain the %following 
% divergence form 
% \begin{align}
% \begin{aligned}
% &\partial_t\Big( \, \frac{1}{r+1} \left|\, \phi \, \right|^{r+1} \, \Big) \\[4pt]
% &+\partial _x \Big( \, 
%   |\, \phi \, |^{r-1} \, \phi \, \bigl( \, f(U+\phi) - f(U) \, \bigr) \, 
%   \Big) \\
% &+\partial _x \Big( \, 
%   - r \, \int _{0}^{\phi} 
%   \bigl( \, f(U+\eta) - f(U) \, \bigr) \, |\, \eta \, |^{r-1} 
%     \, \mathrm{d}\eta 
%   \, \Big) \\
% &+\partial _x \Big( \, 
%   -|\, \phi \, |^{r-1} \, \phi \, 
%   \Big( \, 
%   \sigma\bigl( \, 
%   \partial_x U + \partial_x \phi \, \bigr) - 
%     \sigma\bigl( \, 
%     \partial_x U  \, 
%     \bigr)
%     \, \Big)
%   \, \Big) \\
% &+r\,  \int _{0}^{\phi} 
%   \bigl( \, f'(U+\eta) - f'(U) \, \bigr) \, |\, \eta \, |^{r-1} 
%    \partial_x U \, \mathrm{d}\eta \\ 
% &+r \, |\, \phi \, |^{r-1} \, \partial_x \phi \, 
%   \Big( \, 
%   \sigma\bigl( \, 
%   \partial_x U + \partial_x \phi \, \bigr) - 
%     \sigma\bigl( \, 
%     \partial_x U  \, 
%     \bigr)
%     \, \Big)
% = |\, \phi \, |^{r-1}\, \phi \, 
%   \partial_x \left( \, \sigma\bigl( \, \partial_x U  \, \bigr)  \, \right). 
% \end{aligned}
% \end{align}
% Integrating (4.1) with respect to $x$, we have
\begin{align}
\begin{aligned}
&\frac{1}{r+1} \, 
% \frac{\mathrm{d}}{\mathrm{d}t} \, 
\frac{\mathrm{d}}{\mathrm{d}t} \, 
\Vert \, \phi(t) \, \Vert_{L^{r+1}}^{r+1} 
+r\,  \int
\int _{0}^{\phi} 
  \bigl( \, f'(\eta+U) - f'(U) \, \bigr) \, |\, \eta \, |^{r-1} \partial_x U
   \, \mathrm{d}\eta \, \mathrm{d}x\\
&\quad\quad+r\,  \int ^{\infty }_{-\infty } 
|\, \phi \, |^{r-1} \, \partial_x \phi \, 
  \big( \, 
  \sigma\bigl( \, 
  \partial_x U + \partial_x \phi \, \bigr) - 
    \sigma\bigl( \, 
    \partial_x U  \, 
    \bigr)
    \, \big)
\, \mathrm{d}x\\
&\quad=\int 
 |\, \phi \, |^{r-1}\, \phi \, 
  \partial_x \left( \, \sigma\bigl( \, \partial_x U  \, \bigr)  \, \right)
 \, \mathrm{d}x. 
\end{aligned}
\end{align}
We estimate the right-hand side of (4.1) by the H\"{o}lder's inequality as 
\begin{equation}
\left| \,  
\int
 |\, \phi \, |^{r-1}\, \phi \, 
  \partial_x \left( \, \sigma\bigl( \, \partial_x U  \, \bigr)  \, \right)
 \, \mathrm{d}x
\, \right|
\leq \Vert \, \phi
     \, \Vert_{L^{r+1}}^{r} \, 
     \left|\left| \,    
     \partial_x \left( \, \sigma\bigl( \, \partial_x U  \, \bigr)  \, \right)
     \, \right| \right|_{L^{r+1}}.
\end{equation} 
Note that by the assumptions (1.2),(1.5) on $f$ and $\sigma$,
the second and third terms on the left side of (4.1) are non-negative.
Then, substituting (4.2) into (4.1), we have 
\begin{equation}
 \frac{\mathrm{d}}{\mathrm{d}t} \, 
%\frac{d}{dt} \,
\,\Vert \, \phi
(t) \, \Vert_{L^{r+1}}
\leq \left|\left| \,    
     \partial_x \left( \, \sigma\bigl( \, \partial_x U  \, \bigr)  \, \right)(t) 
     \, \right| \right|_{L^{r+1}}.
\end{equation} 
Integrating (4.3) with respect to $t$, 
we have for any compact set $K\subset \mathbb{R}$, 
\begin{equation}
\Vert \, \phi(t) \, \Vert_{L^{r+1}(K)}
\leq \Vert \, \phi(t) \, \Vert_{L^{r+1}}
\leq \Vert \, \phi_{0} \, \Vert_{L^{r+1}}
     + \int _{0}^{t} \left|\left| \,    
     \partial_x \left( \, \sigma\bigl( \, \partial_x U  \, \bigr)  \, \right)(\tau) 
     \, \right| \right|_{L^{r+1}}\, \mathrm{d} \tau.
\end{equation} 
Taking the limit $r \rightarrow \infty$ in (4.4), we immediately have 
\begin{equation}
\max_{x \in K}
| \, \phi(t, \, x) \, |
\leq \Vert \, \phi_{0} \, \Vert_{L^{\infty}}
     + \int _{0}^{t} \left|\left| \,    
     \left( \, \sigma'\bigl( \, \partial_x U  \, \bigr)  
     \, \partial_x^2 U \, \right)(\tau) 
     \, \right| \right|_{L^{\infty}}\, \mathrm{d} \tau.
\end{equation} 
Because the compact set $K\subset \mathbb{R}$ is arbitrary, 
we obtain 
\begin{equation}
\sup_{x \in \mathbb{R}}
| \, \phi(t, \, x) \, |
\leq %\sup_{x \in \mathbb{R}}
%     | \, \phi_{0} (x) \, |
     \Vert \, \phi_{0} \, \Vert_{L^{\infty}}
     + C \, \int _{0}^{t} \left|\left| \,    
     \, \partial_x^2 U(\tau) 
     \, \right| \right|_{L^{\infty}}\mathrm{d}\tau.
\end{equation} 
Since $\| \partial_x^2 U(\cdot)\|_{L^{\infty}} \in L_t^1(0,\, \infty)$
by Lemma 2.2, the proof of Lemma 4.1 is completed.
% Thus, by using Lemma 2.2, we complete the proof of Lemma 4.1. 

\medskip

Next we show Proposition 4.1. 

\medskip

{\bf Proof of Proposition 4.1.}
Taking $r=1$ in (4.1), and using the assumptions (1.2), (1.3), (1.5)
together with 
Lemma 4.1, we have 
\begin{align}
\begin{aligned}
&\frac{1}{2} \, 
\frac{\mathrm{d}}{\mathrm{d}t} \, 
\Vert \, \phi(t) \, \Vert_{L^{2}}^{2} 
+C^{-1}_0\,  \int
\phi^2\, \partial_x U \, \mathrm{d}x+C^{-1}_0\,  \int 
<\partial_x \phi>^{p-1} 
|\, \partial_x\phi \, |^2 
\, \mathrm{d}x\\
&\quad\leq \bigg| \,  
\int 
 \phi \, 
  \partial_x \left( \, \sigma\bigl( \, \partial_x U  \, \bigr)  \, \right)
 \, \mathrm{d}x
\, \bigg|. 
\end{aligned}
\end{align}
% We estimate the third term on the left-hand side of (4.8) as 
% \begin{align}
% \begin{aligned}
% &\int ^{\infty }_{-\infty } 
% <\, \partial_x \phi \, >^{p-1} 
% |\, \phi \, |^2 
% \, \mathrm{d}x\\
% &\geq C^{-1} \, \Vert \, \partial_x \phi(t) \, 
%        \Vert_{L^{2}(\{ \, |\partial_x \phi|<1 \, \})}^{2} 
%      + C^{-1} \, \Vert \, \partial_x \phi(t) \, 
%        \Vert_{L^{p+1}(\{ \, |\partial_x \phi|\geq1 \, \})}^{p+1}.  
% \end{aligned}
% \end{align}
% We also devide the integral region to the right-hand side of (4.8) as 
% \begin{align}
% \begin{aligned}
% &\left| \,  
% \int ^{\infty}_{-\infty} 
%  \phi \, 
%   \partial_x \left( \, \sigma\bigl( \, \partial_x U  \, \bigr)  \, \right)
%  \, \mathrm{d}x
% \, \right|\\
% &\leq \left( \, \int_{|\partial_x \phi|<1} 
%      + \int_{|\partial_x \phi|\geq1}  \, \right) \, 
%       | \, \phi \, | \big| \, 
%       \partial_x \left( \, \sigma\bigl( \, \partial_x U  \, \bigr)  \, \right)  
%       \, \big| 
%       \, \mathrm{d}x =: I_{1} + I_{2}. 
% \end{aligned}
% \end{align}
% To estimate (4.10), 
% we construct the following Sobolev-type interpolation inequalities as follows 
% (see \cite{yoshida2}, \cite{yoshida3}, \cite{yoshida4}). 
In order to estimate the right hand side of (4.7), we prepare
the following

\medskip

\noindent
{\bf Lemma 4.2.}\quad {\it
For $g \in H^2$, it holds that
\begin{equation}
% \sup_x |\,g(x)\,| 
\|\,g\,\|_{L^{\infty}}^2
\le C\,Q_g^{\frac 1 4}\,\|\,g\,\|_{L^2}^{\frac 1 2}+
C\,Q_g^{\frac {1}{ 3p+1}}\,\|\,g\,\|_{L^2}^{\frac{2p}{3p+1}},
\end{equation}
where 
$$
Q_g := \int <\partial_x g>^{p-1}\,|\,\partial_x g\,|^2\, \mathrm{d}x.
$$
}
\medskip

{\bf Proof of Lemma 4.2.}\quad
We first note that
\begin{equation}
Q_g \sim 
\int_{|\partial_x g|<1}
|\, \partial_x g \, |^2 \, \mathrm{d}x
+
\int_{|\partial_x g|>1}
|\, \partial_x g \, |^{p+1} \, \mathrm{d}x.
\end{equation}
By simple calculation, we have
\begin{equation}
% \sup_x |\,g(x)\,|^2 \le
\|\,g\,\|_{L^{\infty}} \le
\int 2\,|\,g\,|\,|\,\partial_x g\,| \,\mathrm{d}x= \int_{|\,\partial_x g\,|<1} +
\int_{|\,\partial_x g\,|>1} =: I_1+I_2.
\end{equation}
We estimate each $I_i\,(i=1,2)$ by using the H\"{o}lder and Young inequalities as follows:
\begin{equation}
I_1 \le 2 \,\Big(\,\int_{|\,\partial_x g\,|<1}
|\, \partial_x g \, |^2 \, \mathrm{d}x\,\Big)^{\frac 1 2}
\Big(\,\int_{|\,\partial_x g\,|<1}
|\, g \, |^2 \, \mathrm{d}x\,\Big)^{\frac 1 2}
\le C\,Q_g^{\frac 1 2}\,\|\,g\,\|_{L^2};
\end{equation}
\begin{align}
\begin{aligned}
I_2&\leq 
2 \,\Big(\,
\int_{|\,\partial_x g\,|>1}
|\, \partial_x g \, |^{p+1} \, \mathrm{d}x\,
\Big)^{\frac{1}{p+1}}
\Big(\int_{|\,\partial_x g\,|>1}
|\, g \, |^{\frac{p+1}{p}} \, \mathrm{d}x\,\Big)^{\frac{p}{p+1}}
\\
&\leq 2\,Q_g^{\frac{1}{p+1}}\,\|\,g\,\|_{L^2}^{\frac{2p}{p+1}}\,\|\,g\,\|_{L^{\infty}}^{\frac{1-p}{p+1}}
\\
&\leq \epsilon \,\|\,g\,\|_{L^{\infty}}^2+ C_{\epsilon}\,
Q_g^{\frac{2}{3p+1}}\,\|\,g\,\|_{L^2}^{\frac{4p}{3p+1}}\quad (\epsilon>0).          
\end{aligned}
\end{align}
Thus, substituting (4.11) and (4.12) into (4.10), and choosing 
$\epsilon$ suitably small, we complete the proof of Lemma 4.2.

\medskip

Let us turn to the estimate of the right hand side of (4.7) by
using Lemma 2.2 and Lemma 4.2 as
\begin{align}
\begin{aligned}
&\bigg| \,  
\int 
 \phi \, 
  \partial_x \big( \, \sigma ( \, \partial_x U  \,)  \, \big)
 \, \mathrm{d}x
\, \bigg|
\le C\,\|\,\phi\,\|_{L^{\infty}}\|\,\partial_x^2U\,\|_{L^{1}}
\\
& \quad \le
\epsilon \,Q_{\phi} + C_{\epsilon}\,\|\,\phi\,\|_{L^{2}}^{\frac{2}{3}}
\Big(\,\|\,\partial_x^2U\,\|_{L^{1}}^{\frac{4}{3}}
+\|\,\partial_x^2U\,\|_{L^{1}}^{\frac{3p+1}{3p}}\,\Big)
\\
& \quad \le
\epsilon \,Q_{\phi} + C_{\epsilon}\,(\,1+\|\,\phi\,\|_{L^{2}}^{2}\,)\,(1+t)^{-\frac{4}{3}}
\qquad (\epsilon>0).
\end{aligned}
\end{align}

Substituting (4.13) into (4.7), 
choosing $\epsilon$ suitably small, 
and using the Gronwall inequality, 
we obtain the desired {\it a priori} estimate for $\phi$. 
Thus, the proof of Proposition 4.1 is completed. 

\bigskip 

\noindent
%%%%%%%%%%%%%%%%%%%%%%%%%%%%%%%%%%%%%%%%%%%%%%%%%%%%%%%%%%%%%%%%%%%%%%%%%%%%%
\section{A priori estimates I\hspace{-.1em}I}
In this section, 
we proceed to the {\it a priori} estimate 
for the derivative $\partial_x \phi$. 

\medskip

\noindent
{\bf Proposition 5.1.}\quad {\it
For $0<p<1$, there exists a positive constant $C_0$
such that 
\begin{align*}
\begin{aligned}
& \int 
     <\partial_x \phi>^{p-1} 
|\, \partial_x \phi \, |^2 
\, \mathrm{d}x
   + \int^{t}_{0} \int
     <\partial_x \phi>^{2\,(p-1)} 
|\, \partial_x^2 \phi \, |^2 
\, \mathrm{d}x\mathrm{d}\tau \\[5pt]
% &+ \int^{t}_{0} \int ^{\infty }_{-\infty } 
%   \Big( \, 
%   \sigma'\bigl( \, 
%   \partial_x U + \partial_x \phi \, \bigr) - 
%     \sigma'\bigl( \, 
%     \partial_x U  \, 
%     \bigr)  \, 
%     \Big)^2 \, 
%     ( \, \partial_x^2 U \, )^2  
%     \, \mathrm{d}x\mathrm{d}\tau \\
%& \qquad \qquad \qquad \qquad \qquad \qquad 
% &\leq C_0 ( \, \phi_{0}, \: \partial_x\phi_{0} \, )\, 
%      \Bigl( \, 1+
%      \ll \, \partial_x\phi \, \gg _{L^{\infty}_{t, \, x}}^{1-p}
%      \, \Bigr) 
% \quad \big( \, t \in [\, 0, \: T\, ] \, \big), 
&\leq C_0 
     \ll \partial_x\phi \gg _{\infty}^{1-p}
\qquad \big( \, t \in [\, 0, \: T\, ] \,\big), 
\end{aligned}
\end{align*}
where 
$
% \displaystyle{
% \ll \, v \, \gg _{L^{\infty}_{t, \, x}}
% := \esssup\displaylimits_{t\geq0, \; x \in \mathbb{R}} \, <\, v(t, \,x) \,>}. 
\displaystyle{
\ll v  \gg _{\infty}
:= \esssup\displaylimits_{t\in [0,T], \; x \in \mathbb{R}} \, < v(t, \,x) >}. 
$
}

\medskip

{\bf Proof of Proposition 5.1.}
Multiplying the equation in (3.2) by 
$$
- \partial_x 
  \big( \, 
  \sigma( \, 
  \partial_x U + \partial_x \phi \,) - 
    \sigma ( \, 
    \partial_x U  \,)
    \, \big),
$$
% we obtain the %following 
% divergence form 
% \begin{align}
% \begin{aligned}
% &\partial_t\left(\, 
% \int _{0}^{\partial_x \phi} 
% \Big( \, 
%   \sigma\bigl( \, 
%   \partial_x U + \eta \, \bigr) - 
%     \sigma\bigl( \, 
%     \partial_x U  \, 
%     \bigr)  \, 
%     \Big) \, \mathrm{d}\eta 
%  \, \right) \\
% &+\partial _x \bigg( \, -\partial_t \phi \, \Big( \, 
%   \sigma\bigl( \, 
%   \partial_x U + \partial_x \phi \, \bigr) - 
%     \sigma\bigl( \, 
%     \partial_x U  \, 
%     \bigr)  \, 
%     \Big) \, \bigg)\\
% &-\partial _x \bigl( \, f(U+\phi) - f(U) \, \bigr) \, 
%   \partial _x \Big( \, 
%   \sigma\bigl( \, 
%   \partial_x U + \partial_x \phi \, \bigr) - 
%     \sigma\bigl( \, 
%     \partial_x U  \, 
%     \bigr)  \, 
%     \Big)\\
% &-\int _{0}^{\partial_x \phi} 
%   \Big( \, 
%   \sigma'\bigl( \, 
%   \partial_x U + \eta \, \bigr) - 
%     \sigma'\bigl( \, 
%     \partial_x U  \, 
%     \bigr)  \, 
%     \Big) \, \partial_t \partial_x U \, \mathrm{d}\eta  \\
% &+\bigg( \, \partial _x \Big( \, 
%   \sigma\bigl( \, 
%   \partial_x U + \partial_x \phi \, \bigr) - 
%     \sigma\bigl( \, 
%     \partial_x U  \, 
%     \bigr)  \, 
%     \Big)
%     \, \bigg)^2  \\
% &= -\partial _x \Big( \, 
%   \sigma\bigl( \, 
%   \partial_x U + \partial_x \phi \, \bigr) - 
%     \sigma\bigl( \, 
%     \partial_x U  \, 
%     \bigr)  \, 
%     \Big) \, 
%   \partial_x \left( \, \sigma\bigl( \, \partial_x U  \, \bigr)  \, \right). 
% \end{aligned}
% \end{align}
and integrating the resultant formula with respect to $x$, we have,
after integration by parts,
\begin{align}
\begin{aligned}
&\frac{\mathrm{d}}{\mathrm{d}t} \, 
\int 
\int _{0}^{\partial_x \phi} 
\big( \, 
  \sigma ( \, 
  \partial_x U + \eta \,) - 
    \sigma ( \, 
    \partial_x U  \,)  \, 
    \big) \, \mathrm{d}\eta\mathrm{d}x\\
&-\int 
  \partial _x \bigl( \, f(U+\phi) - f(U) \, \bigr) \, 
  \partial _x \big( \, 
  \sigma ( \, 
  \partial_x U + \partial_x \phi \,) - 
    \sigma ( \, 
    \partial_x U  \,)  \, 
    \big) 
    \, \mathrm{d}x\\
&-\int  
  \big( \, 
  \sigma ( \, \partial_x U + \partial_x\phi \, ) - 
  \sigma ( \, \partial_x U \, ) -
    \sigma'( \, \partial_x U  \,)\partial_x\phi  \, 
    \big) \partial_t \partial_x U  \, \mathrm{d}x\\
&+\int 
  \big| \, \partial _x \big( \, 
  \sigma ( \, 
  \partial_x U + \partial_x \phi \,) - 
    \sigma ( \, 
    \partial_x U  \,)  \, 
    \big)\, \big|^2  
    \, \mathrm{d}x\\
&= -\int 
  \partial _x \big( \, 
  \sigma ( \, 
  \partial_x U + \partial_x \phi \,) - 
    \sigma ( \, 
    \partial_x U  \,)  \, 
    \big) \, 
  \partial_x \big( \, \sigma ( \, \partial_x U  \,)  \, \big)
\, \mathrm{d}x. 
\end{aligned}
\end{align}
By using the Young inequality, we estimate the second term 
on the left-hand side of (5.1) as 
\begin{align}
\begin{aligned}
& \left| \,  
  \int 
  \partial _x \big( \, f(U+\phi) - f(U) \, \big) \, 
  \partial _x \big( \, 
  \sigma ( \, 
  \partial_x U + \partial_x \phi \,) - 
    \sigma ( \, 
    \partial_x U  \,)  \, 
    \big) 
    \, \mathrm{d}x
\, \right|\\
& \leq \epsilon \, 
       \int 
  \big| \, \partial _x \big( \, 
  \sigma ( \, 
  \partial_x U + \partial_x \phi \,) - 
    \sigma ( \, 
    \partial_x U  \,)  \, 
    \big)
    \, \big|^2  
    \, \mathrm{d}x\\
& \qquad 
    + C_{\epsilon} \, 
      \int 
       \big| \, \partial _x ( \, f(U+\phi) - f(U) \,\big)
    \, \big|^2  
    \, \mathrm{d}x\qquad (\epsilon>0). 
\end{aligned}
\end{align}
Similarly, the right-hand side of (5.1) is estimated as 
\begin{align}
\begin{aligned}
& \left| \,
\int 
  \partial _x \big( \, 
  \sigma ( \, 
  \partial_x U + \partial_x \phi \,) - 
    \sigma ( \, 
    \partial_x U  \,)  \, 
    \big) \, 
  \partial_x \big( \, \sigma ( \, \partial_x U  \,)  \, \big)
\, \mathrm{d}x 
%   \int ^{\infty }_{-\infty } 
%   \Big( \, \partial _x \Big( \, 
%   \sigma\bigl( \, 
%   \partial_x U + \partial_x \phi \, \bigr) - 
%     \sigma\bigl( \, 
%     \partial_x U  \, 
%     \bigr)  \, 
%     \Big) \, 
%   \partial_x \left( \, \sigma\bigl( \, \partial_x U  \, \bigr)  \, \right)
% \, \mathrm{d}x
\, \right|\\
& \leq \epsilon \, 
       \int
  \big| \, \partial _x \big( \, 
  \sigma ( \, 
  \partial_x U + \partial_x \phi \,) - 
    \sigma ( \, 
    \partial_x U  \,)  \, 
    \big)
    \, \big|^2  
    \, \mathrm{d}x\\
& \qquad 
    + C_{\epsilon} \, 
      \int 
%   \big| \, 
%        \partial_x \big( \, \sigma ( \, \partial_x U  \,)  \, \big)\, \big|^2  
  | \, \partial_x^2 U  \,|^2  
    \, \mathrm{d}x \qquad (\epsilon >0). 
\end{aligned}
\end{align}
% Noting by using the Tayler formula that 
The third term
on the left side of (5.1) is estimated 
by the Tayler formula, the uniform boundedness of 
$\sigma'$ for $0<p<1$, and Lemma 2.1 as
\begin{align}
\begin{aligned}
&\bigg|  \, \int \big( \, 
  \sigma ( \, \partial_x U + \partial_x\phi \, ) - 
  \sigma ( \, \partial_x U \, ) -
    \sigma'( \, \partial_x U  \,) \, \partial_x\phi  \, 
    \big) \,\partial_t \partial_x U\,\mathrm{d}x  \, \bigg|
\\[5pt]
&\le \int \big|\, \big(\,\sigma'( \, 
  \partial_x U + \theta \, \partial_x \phi \,)
  - \sigma'( \, \partial_x U  \,)\,\big)\big|\,|\,\partial_x\phi\,|\,
|\,\partial_t \partial_x U\,| \,\mathrm{d}x
\\[5pt]
&\le C\,\int | \, \partial_x \phi \, |^2\, \mathrm{d}x
 \qquad ( \,  \exists \theta=\theta(t, \, x)\in (0, \,1) \,).
\end{aligned}
\end{align}
% \begin{align}
% \begin{aligned}
% &\int ^{\infty }_{-\infty } \int _{0}^{\partial_x \phi} 
%   \Big( \, 
%   \sigma'\bigl( \, 
%   \partial_x U + \eta \, \bigr) - 
%     \sigma'\bigl( \, 
%     \partial_x U  \, 
%     \bigr)  \, 
%     \Big) \, \mathrm{d}\eta \, ( \, \partial_t \partial_x U \, ) \, \mathrm{d}x\\
% &= \Big( \, 
%   \sigma\bigl( \, 
%   \partial_x U + \partial_x \phi \, \bigr) - 
%     \sigma\bigl( \, 
%     \partial_x U  \, 
%     \bigr)  - 
%     \sigma'\bigl( \, 
%     \partial_x U  \, 
%     \bigr) \, 
%     \partial_x \phi
%     \, \Big) \, 
%     \partial_t \partial_x U\\
% &= \frac{1}{2} \, \sigma''\bigl( \, 
%   \partial_x U + \theta \, \partial_x \phi \, \bigr) \, 
%   ( \, \partial_x \phi \, )^2 \, \partial_t \partial_x U 
%   \qquad \big( \,  \exists \theta=\theta(t, \, x) \, \big),  
% \end{aligned}
% \end{align}
% the third term on the left-hand side of (5.2) is estimated as 
% \begin{align}
% \begin{aligned}
% & \left| \,  
%   \int ^{\infty }_{-\infty } \int _{0}^{\partial_x \phi} 
%   \Big( \, 
%   \sigma'\bigl( \, 
%   \partial_x U + \eta \, \bigr) - 
%     \sigma'\bigl( \, 
%     \partial_x U  \, 
%     \bigr)  \, 
%     \Big) \, \mathrm{d}\eta \, ( \, \partial_t \partial_x U \, ) \, \mathrm{d}x
% \, \right|\\
% & \leq C \, 
%        \| \, \partial_t \partial_x U(t) \, \|_{L^{\infty}} \, 
%        \| \, \partial_x \phi(t) \, \|_{L^2}^2. 
% \end{aligned}
% \end{align}
Substituting (5.2), (5.3), and (5.4) into (5.1), and choosing 
$\epsilon$ suitably small, we have 
\begin{align}
\begin{aligned}
&\frac{\mathrm{d}}{\mathrm{d}t} \, 
\int  
\int _{0}^{\partial_x \phi} 
\big( \, 
  \sigma ( \, 
  \partial_x U + \eta \, ) - 
    \sigma ( \, 
    \partial_x U  \,)  \, 
    \big) \, \mathrm{d}\eta\mathrm{d}x\\ 
&+\int 
  \big| \, \partial _x \big( \, 
  \sigma ( \, 
  \partial_x U + \partial_x \phi \,) - 
    \sigma ( \, 
    \partial_x U  \,)  \, 
    \big)\big|^2  
    \, \mathrm{d}x\\
&\leq C \, \int 
       \big| \, \partial _x \big( \, f(U+\phi) - f(U) \, \big)
    \, \big|^2  
    \, \mathrm{d}x \\
& \quad 
%     + C \, \left\| \,  
%        \partial_x \left( \, \sigma\bigl( \, \partial_x U  \, \bigr)\, \right)(t) 
%        \, \right\|_{L^2}^2
%     + C \, \| \, \partial_x \phi(t) \, \|_{L^2}^2.
% +C \int \big(
%      \big| \,  
%        \partial_x \big( \, \sigma ( \, \partial_x U  \,)\, \big|^2
%     +  | \, \partial_x \phi \,|^2
+C \, \int \big(\,
     | \,  
       \partial_x^2 U \,|^2
    +  | \, \partial_x \phi \,|^2
\,\big)\, \mathrm{d}x.
\end{aligned}
\end{align}
% provided $\epsilon$ is suitably small. 
% By using the uniform boundedness of $\phi$ (Lemma 4.1) and $U$ (Lemma 2.2), 
By Lemma 2.2 and Lemma 4.1,
we estimate the first term on the right-hand side of (5.5) as 
% \begin{align}
% \begin{aligned}
% &\int 
%        \big| \, \partial _x \big( \, f(U+\phi) - f(U) \, )
%     \, \big|^2  
%     \, \mathrm{d}x\\
% &\leq C \, \| \, \partial_x U(t) \, \|_{L^{\infty}} \, \int ^{\infty }_{-\infty } 
%   \phi^2 \, 
%   \partial_x U
%   \, \mathrm{d}x + C \, \| \, \partial_x \phi(t) \, \|_{L^2}^2. 
% \end{aligned}
% \end{align}
\begin{equation}
\int 
       \big| \, \partial _x \big( \, f(U+\phi) - f(U) \, )
    \, \big|^2  
    \, \mathrm{d}x
\leq C \, \int \big(\,
  \phi^2 \, 
  \partial_x U  + | \, \partial_x \phi \,|^2 \big)\, \mathrm{d}x . 
\end{equation}
Similarly, the second term on the left-hand side of (5.5) 
is estimated as
\begin{align}
\begin{aligned}
&\int  
  \big|\, \partial _x \big( \, 
  \sigma ( \, 
  \partial_x U + \partial_x \phi \,) - 
    \sigma ( \, 
    \partial_x U  \, 
    )  \, 
    \big)\big|^2  
    \, \mathrm{d}x\\
&\ge \int 
  \big| \, 
  \sigma'\big( \, 
  \partial_x U + \partial_x \phi \, \big)  \, 
    \big|^2 \, 
    |\, \partial_x^2 \phi \, |^2  
    \, \mathrm{d}x\\
& \quad
  - \int 
  \big| \, 
  \sigma'\big( \, 
  \partial_x U + \partial_x \phi \, \big) - 
    \sigma'\big( \, 
    \partial_x U  \, 
    \big)  \, 
    \big|^2 \, 
    |\, \partial_x^2 U \,|^2  
    \, \mathrm{d}x\\
&\ge \int 
  \big| \, 
  \sigma'\big( \, 
  \partial_x U + \partial_x \phi \, \big)  \, 
    \big|^2 \, 
    |\, \partial_x^2 \phi \, |^2  
    \, \mathrm{d}x
  -  C \, \int 
    |\, \partial_x^2 U \,|^2  
    \, \mathrm{d}x.
% &
% + 2 \int ^{\infty }_{-\infty } 
%   \sigma'\bigl( \, 
%   \partial_x U + \partial_x \phi \, \bigr)  \, 
%   \partial_x^2 \phi \, 
%   \Big( \, 
%   \sigma'\bigl( \, 
%   \partial_x U + \partial_x \phi \, \bigr) - 
%     \sigma'\bigl( \, 
%     \partial_x U  \, 
%     \bigr)  \, 
%     \Big)\, \partial_x^2 U 
%     \, \mathrm{d}x. 
\end{aligned}
\end{align}
% By using the Young inequality and Lagrange's mean-value theorem, 
% the third term on the right-hand side of (5.9) is estimated as 
% \begin{align}
% \begin{aligned}
% & \left| \,  
%   2 \int ^{\infty }_{-\infty } 
%   \sigma'\bigl( \, 
%   \partial_x U + \partial_x \phi \, \bigr)  \, 
%   \partial_x^2 \phi \, 
%   \Big( \, 
%   \sigma'\bigl( \, 
%   \partial_x U + \partial_x \phi \, \bigr) - 
%     \sigma'\bigl( \, 
%     \partial_x U  \, 
%     \bigr)  \, 
%     \Big)\, \partial_x^2 U 
%     \, \mathrm{d}x
% \, \right|\\
% & \leq \epsilon \, 
%        \int ^{\infty }_{-\infty } 
%   \Big( \, 
%   \sigma'\bigl( \, 
%   \partial_x U + \partial_x \phi \, \bigr)  \, 
%     \Big)^2 \, 
%     ( \, \partial_x^2 \phi \, )^2  
%     \, \mathrm{d}x
%     + C_{\epsilon} \, 
%       \int ^{\infty }_{-\infty } 
%       ( \, \partial_x \phi \, )^2
%       ( \, \partial_x^2 U \, )^2
%     \, \mathrm{d}x. 
% \end{aligned}
% \end{align}
% Noting 
% \begin{equation}
% \int ^{\infty }_{-\infty } 
%   \phi^2 \, 
%   \partial_x U
%   \, \mathrm{d}x\in L_{t}^1(0, \: \infty), \quad 
% \left\| \,  
%        \partial_x \left( \, \sigma\bigl( \, \partial_x U  \, \bigr)\, \right)
%        \, \right\|_{L^2}^2 \in L_{t}^1(0, \: \infty), 
% \end{equation}
% from Proposition 4.1 and Lemma 2.2, and 
Furthermore, by the assumptions (1.2),(1.3), it holds
\begin{equation}
\int 
  | \, 
  \sigma' ( \, 
  \partial_x U + \partial_x \phi \, )  \,|^2 \, 
    | \, \partial_x^2 \phi \, |^2  
    \, \mathrm{d}x
\sim 
     \int 
     <\partial_x \phi>^{2\,(p-1)} 
     |\, \partial_x^2 \phi \, |^2 
     \, \mathrm{d}x, 
\end{equation}
and
\begin{equation}
\int 
\int _{0}^{\partial_x \phi} 
\big( \, 
  \sigma ( \, 
  \partial_x U + \eta \,) - 
    \sigma ( \, 
    \partial_x U  \, 
    )  \, 
    \big) \, \mathrm{d}\eta\mathrm{d}x
\sim 
     \int 
     <\partial_x \phi>^{p-1} 
     |\, \partial_x \phi \, |^2 
     \, \mathrm{d}x. 
\end{equation}
Hence, substituting (5.6) and (5.7) into (5.5),
and integrating with respect to $t$, 
we have 
% substituting (5.8), (5.9) and (5.10) into (5.7), 
% choosing $\epsilon$ suitably small, 
% and integrating with respect to $t$, 
% we have 
\begin{align}
\begin{aligned}
& \int 
     <\partial_x \phi>^{p-1} 
|\, \partial_x \phi \, |^2 
\, \mathrm{d}x
   + \int^{t}_{0} \int 
     <\partial_x \phi>^{2\,(p-1)} 
|\, \partial_x^2 \phi \, |^2 
\, \mathrm{d}x\mathrm{d}\tau \\
% &+ \int^{t}_{0} \int ^{\infty }_{-\infty } 
%   \Big( \, 
%   \sigma'\bigl( \, 
%   \partial_x U + \partial_x \phi \, \bigr) - 
%     \sigma'\bigl( \, 
%     \partial_x U  \, 
%     \bigr)  \, 
%     \Big)^2 \, 
%     ( \, \partial_x^2 U \, )^2  
%     \, \mathrm{d}x\mathrm{d}\tau \\
%& \qquad \qquad \qquad \qquad \qquad \qquad 
% &\leq C_0 
%      \Big( \, 1+
%      \int^{t}_{0} 
%      \| \, \partial_x \phi(\tau) \, \|_{L^2}^2 \, \mathrm{d}\tau
%      \, \Big).  
&\leq C_0 \,
     \Big( \, 1+
     \int^{t}_{0} 
     \int \big(\,
  \phi^2 \, 
  \partial_x U  + | \, \partial_x \phi \,|^2 + |\,\partial_x^2 U\,|^2
\big)\, \mathrm{d}x
\mathrm{d}\tau
     \, \Big).
\end{aligned}
\end{align}
Finally, by Lemma 2.1 and Proposition 4.1, it holds that
\begin{equation}
     \int^{t}_{0} 
     \int \big(\,
  \phi^2 \, 
  \partial_x U  + |\,\partial_x^2 U\,|^2
\big)\, \mathrm{d}x
\mathrm{d}\tau \le C_0,
\end{equation}
and 
% Noting from Proposition 4.1 that 
\begin{align}
\begin{aligned}
 \int^{t}_{0} \int
     | \, \partial_x \phi(\tau) \,|^2 \, \mathrm{d}x\mathrm{d}\tau
 &\leq  
      \int^{t}_{0} \int \ll\partial_x\phi \gg _{\infty}^{1-p}
      <\partial_x \phi>^{p-1} 
      |\, \partial_x \phi \, |^2 
      \, \mathrm{d}x\mathrm{d}\tau 
\\
&\leq C_0 
      \ll \partial_x\phi \gg _{\infty}^{1-p}.
\end{aligned}
\end{align}
Substituting (5.11) and (5.12) into 
(5.10), we obtain the desired {\it a priori} estimate for 
$\partial_x \phi$. 
Thus, the proof of Proposition 5.1 is completed.

\bigskip 

\noindent
%%%%%%%%%%%%%%%%%%%%%%%%%%%%%%%%%%%%%%%%%%%%%%%%%%%%%%%%%%%%%%%%%%%%%%%%%%%%%
\section{A priori estimates I\hspace{-.1em}I\hspace{-.1em}I}

In this section, 
we further show the {\it a priori} estimate 
for $\partial_x^2 \phi$, establish the uniform
boundedness of $\|\,\partial_x \phi\,\|_{L^{\infty}}$, 
and then accomplish the proof of Theorem 3.3. 

\medskip

\noindent
% {\bf Proposition 6.1.}\quad {\it
% For any $3/7<p<1$, there exists a positive constant $C$
% such that 
% \begin{align*}
% \begin{aligned}
% & \int ^{\infty }_{-\infty } 
%      <\, \partial_x \phi \, >^{2\,(p-1)} 
% |\, \partial_x^2 \phi \, |^2 
% \, \mathrm{d}x
%    + \int ^{\infty }_{-\infty } 
%   \Big( \, 
%   \sigma'\bigl( \, 
%   \partial_x U + \partial_x \phi \, \bigr) - 
%     \sigma'\bigl( \, 
%     \partial_x U  \, 
%     \bigr)  \, 
%     \Big)^2 \, 
%     ( \, \partial_x^2 U \, )^2  
%     \, \mathrm{d}x\\
%    &+ \int^{t}_{0} \int ^{\infty }_{-\infty } 
%      <\, \partial_x \phi \, >^{p-1} 
% |\, \partial_t \partial_x \phi \, |^2 
% \, \mathrm{d}x\mathrm{d}\tau 
% \leq C ( \, \phi_{0}, \: \partial_x\phi_{0}, \: \partial_x^2\phi_{0} \, )
% \quad \big( \, t \in [\, 0, \: T\, ] \, \big). 
% \end{aligned}
% \end{align*}
% }
{\bf Proposition 6.1.}\quad {\it
For $0<p<1$, there exists a positive constant $C_0$
such that 
\begin{align*}
\begin{aligned}
\int 
     <\partial_x \phi >^{2\,(p-1)} 
|\, \partial_x^2 \phi \, |^2 
\, \mathrm{d}x
&+ \int^{t}_{0} \int 
     <\partial_x \phi>^{p-1} 
|\, \partial_t \partial_x \phi \, |^2 
\, \mathrm{d}x\mathrm{d}\tau\\[5pt]
&\leq C_0\ll \partial_x\phi \gg _{\infty}^{3\,(1-p)}
\quad \big( \, t \in [\, 0, \: T\, ] \, \big). 
\end{aligned}
\end{align*}
}

\medskip

{\bf Proof of Proposition 6.1.}\quad
% Differentiating the equation in (3.3) with respect to $x$ 
% and multiplying it by 
Differentiating the equation in (3.2) once with respect to $x$ gives
\begin{align*}
\begin{aligned}
  \partial_t\partial_x\phi
  &+ \partial_x^2 \big( \, f(U+\phi) - f(U) \, \big) 
\\[5pt]
%\quad \qquad 
&- \partial_x^2 
    \big( \, 
    \sigma ( \, \partial_x U + \partial_x \phi \, ) 
    - \sigma ( \, \partial_x U  \, )  \, 
    \big)
    = 
    \partial_x^2 \big( \, \sigma ( \, \partial_x U  \,)  \, \big).
\end{aligned}
\end{align*}
Multiplying (6.1) by
$$
 \partial_t 
  \big( \, 
  \sigma ( \, 
  \partial_x U + \partial_x \phi \,) - 
    \sigma( \, 
    \partial_x U  \,)
    \, \big),
$$
and integrating the resultant formula with respect to $x$, we
obtain, after integration by parts,
% we obtain the %following 
% divergence form 
% \begin{align}
% \begin{aligned}
% &\partial_t\left(\, 
% \frac{1}{2} \, 
% \bigg( \, \partial _x \Big( \, 
%   \sigma\bigl( \, 
%   \partial_x U + \partial_x \phi \, \bigr) - 
%     \sigma\bigl( \, 
%     \partial_x U  \, 
%     \bigr)  \, 
%     \Big)
%     \, \bigg)^2
%  \, \right) \\
% &+\partial _x \bigg( \, -\partial_t \Big( \, 
%   \sigma\bigl( \, 
%   \partial_x U + \partial_x \phi \, \bigr) - 
%     \sigma\bigl( \, 
%     \partial_x U  \, 
%     \bigr)  \, 
%     \Big) \, 
%     \partial _x \Big( \, 
%   \sigma\bigl( \, 
%   \partial_x U + \partial_x \phi \, \bigr) - 
%     \sigma\bigl( \, 
%     \partial_x U  \, 
%     \bigr)  \, 
%     \Big)
%     \bigg)\\
% &+\partial _t \partial _x \phi \, 
%   \partial _t \Big( \, 
%   \sigma\bigl( \, 
%   \partial_x U + \partial_x \phi \, \bigr) - 
%     \sigma\bigl( \, 
%     \partial_x U  \, 
%     \bigr)  \, 
%     \Big)\\
% &+\partial _x^2 \bigl( \, f(U+\phi) - f(U) \, \bigr) \, 
%   \partial _t \Big( \, 
%   \sigma\bigl( \, 
%   \partial_x U + \partial_x \phi \, \bigr) - 
%     \sigma\bigl( \, 
%     \partial_x U  \, 
%     \bigr)  \, 
%     \Big)\\
% &= \partial _t \Big( \, 
%   \sigma\bigl( \, 
%   \partial_x U + \partial_x \phi \, \bigr) - 
%     \sigma\bigl( \, 
%     \partial_x U  \, 
%     \bigr)  \, 
%     \Big) \, 
%   \partial_x^2 \left( \, \sigma\bigl( \, \partial_x U  \, \bigr)  \, \right). 
% \end{aligned}
% \end{align}
% Integrating (6.1) with respect to $x$, we have
\begin{align}
\begin{aligned}
&\frac{1}{2} \, \frac{\mathrm{d}}{\mathrm{d}t} \, 
\int 
\big| \, \partial _x \big( \, 
  \sigma( \, 
  \partial_x U + \partial_x \phi \,) - 
    \sigma( \, 
    \partial_x U  \,)  \, 
    \big)
    \, \big|^2 \, \mathrm{d}x\\ 
&+\int  
  \partial _t \partial _x \phi \, 
  \partial _t \big( \, 
  \sigma( \, 
  \partial_x U + \partial_x \phi \,) - 
    \sigma( \, 
    \partial_x U  \,)  \, 
    \big)
    \, \mathrm{d}x\\
&+\int 
\partial _x^2 \big( \, f(U+\phi) - f(U) \, \big) \, 
  \partial _t \big( \, 
  \sigma( \, 
  \partial_x U + \partial_x \phi \,) - 
    \sigma( \, 
    \partial_x U  \,)  \, 
    \big) \, \mathrm{d}x\\
&= \int 
  \partial _t \big( \, 
  \sigma( \, 
  \partial_x U + \partial_x \phi \,) - 
    \sigma( \, 
    \partial_x U  \,)  \, 
    \big) \, 
  \partial_x^2 \big( \, \sigma( \, \partial_x U  \,)  \, \big)
\, \mathrm{d}x. 
\end{aligned}
\end{align}
First, we estimate the second term on the left-hand side of (6.2)
\begin{align}
\begin{aligned}
&\int  
  \partial _t \partial _x \phi \, 
  \partial _t \big( \, 
  \sigma( \, 
  \partial_x U + \partial_x \phi \,) - 
    \sigma( \, 
    \partial_x U  \, 
    )  \, 
    \big)
    \, \mathrm{d}x\\
&= \int 
  \sigma'\bigl( \, 
  \partial_x U + \partial_x \phi \, \bigr) \, 
  | \, \partial _t \partial _x \phi \, |^2 
  \, \mathrm{d}x\\
&\quad + \int 
   \partial _t \partial _x \phi \, 
   \big( \, 
  \sigma'( \, 
  \partial_x U + \partial_x \phi \,) - 
    \sigma'( \, 
    \partial_x U  \, 
    )  \, 
    \big) \, 
    \partial _t \partial _x U 
    \, \mathrm{d}x,
\end{aligned}
\end{align}
% Also the first term on the right-hand side of (6.3) becomes 
where note that
\begin{equation}
\int 
  \sigma'( \, 
  \partial_x U + \partial_x \phi \,) \, 
  | \, \partial _t \partial _x \phi \, |^2 
  \, \mathrm{d}x
\sim 
     \int 
     <\partial_x \phi >^{p-1} 
     |\, \partial _t \partial_x \phi \, |^2 
     \, \mathrm{d}x. 
\end{equation}
The second term on the right-hand side of (6.3) is
estimated by the Young inequality as
\begin{align}
\begin{aligned}
&\left| \,  
 \int 
   \partial _t \partial _x \phi \, 
   \big( \, 
  \sigma'( \, 
  \partial_x U + \partial_x \phi \,) - 
    \sigma'( \, 
    \partial_x U  \, 
    )  \, 
    \big) \, 
    \partial _t \partial _x U 
    \, \mathrm{d}x
 \, \right|\\
& \leq C \, \int 
       |\, \partial_x \phi \, |\,
        |\, \partial _t \partial_x \phi \, |\,
       |\, \partial _t \partial_x U \, |     
       \, \mathrm{d}x\\
& \leq \epsilon \, \int 
       <\partial_x \phi>^{p-1} 
     |\, \partial _t \partial_x \phi \, |^2 
     \, \mathrm{d}x\\
& \quad 
     + C_{\epsilon} \, \int  
       <\partial_x \phi>^{1-p} 
     |\, \partial_x \phi \, |^2 \, 
     |\, \partial _t \partial_x U \, |^2
     \, \mathrm{d}x\\
& \leq \epsilon \, \int 
       <\partial_x \phi>^{p-1} 
     |\, \partial _t \partial_x \phi \, |^2 
     \, \mathrm{d}x\\
& \quad 
     + C_{\epsilon} \, 
%      \| \, \partial _t \partial_x U(t) \, \|_{L^{\infty}}^2\, 
      \ll \partial_x\phi \gg _{\infty}^{2\,(1-p)} \, %\\
     \int 
       <\, \partial_x \phi \, >^{p-1} 
     |\, \partial_x \phi \, |^2 
     \, \mathrm{d}x.
\end{aligned}
\end{align}
Next, we estimate the third term on the left-hand side of (6.2).
Noting 
\begin{align}
\begin{aligned}
& \left|\, \partial _x^2 \big( \, f(U+\phi) - f(U) \, \big) \, \right| \\[5pt]
& \leq C \, \left( \, 
       |\, \partial_x \phi \, |^2 + |\, \partial_x \phi \,|\, |\, \partial_x U \,|
       + |\, \partial_x^2 \phi \, |
       + |\, \phi \, | \, |\, \partial_x U \, |^2
       + |\, \phi \, | \, |\, \partial_x^2 U \, |
       \, \right),
\end{aligned}
\end{align}
and
% \begin{equation}
% \left|\, \partial _t \big( \, 
%   \sigma( \, 
%   \partial_x U + \partial_x \phi \,) - 
%     \sigma ( \, 
%     \partial_x U  \, 
%     )  \, 
%     \big) \, \right| 
%  \leq C (<\partial_x \phi>^{p-1} \, 
%      |\, \partial _t \partial_x \phi \, |
%        + |\, \partial_x \phi \, | \, 
%      |\, \partial _t \partial_x U \, |),
% \end{equation}
\begin{align}
\begin{aligned}
&\left|\, \partial _t \big( \, 
  \sigma( \, 
  \partial_x U + \partial_x \phi \,) - 
    \sigma ( \, 
    \partial_x U  \, 
    )  \, 
    \big) \, \right| \\[5pt]
& \qquad \qquad \leq C \, <\partial_x \phi>^{p-1} \, 
     |\, \partial _t \partial_x \phi \, |
       + C \, |\, \partial_x \phi \, | \, 
     |\, \partial _t \partial_x U \, |,
\end{aligned}
\end{align}
% we estimate the third term on the right-hand side of (6.2) as 
we have
\begin{align}
\begin{aligned}
&\left| \,  
\int 
\partial _x^2 \big( \, f(U+\phi) - f(U) \, \big) \, 
  \partial _t \big( \, 
  \sigma ( \, 
  \partial_x U + \partial_x \phi \,) - 
    \sigma ( \, 
    \partial_x U  \, 
    )  \, 
    \big) \, \mathrm{d}x
\, \right|\\
&\ \  \leq C \int 
       <\partial_x \phi>^{p-1} 
     |\, \partial _t \partial_x \phi \, | \\
&\quad \quad \times 
     \Big( \, 
       |\, \partial_x \phi \, |^2 + |\, \partial_x \phi \, |\,| \, \partial_x U\,|
       + |\, \partial_x^2 \phi \, |
       + |\, \phi \, |\, \big( \, |\, \partial_x U \,|^2
       + |\, \partial_x^2 U \, |
       \, \big)\,  \Big)
     \, \mathrm{d}x\\[5pt]
&\ \  + C \int 
       |\, \partial_x \phi \, | \, 
     |\, \partial _t \partial_x U \, |\\
&\quad\quad \times 
     \Big( \, 
       |\, \partial_x \phi \, |^2 + |\, \partial_x \phi \, |\,| \, \partial_x U\,|
       + |\, \partial_x^2 \phi \, |
       + |\, \phi \, |\, \big(\, |\, \partial_x U \,|^2
       + |\, \partial_x^2 U \, |
       \, \big) \, \Big)
     \, \mathrm{d}x\\
& =: I_{1}+I_{2}. 
\end{aligned}
\end{align}
Let us estimate each $I_{i}\ (i=1,2)$. 
By using the Young inequality, we have 
\begin{align}
\begin{aligned}
 &I_{1}\leq \epsilon \, \int 
       <\partial_x \phi>^{p-1} 
     |\, \partial _t \partial_x \phi \, |^2 
     \, \mathrm{d}x%\\
%& 
+ C_{\epsilon} \, \int 
    <\partial_x \phi>^{p-1} \\
& \quad \times 
    \Big( \, 
       |\, \partial_x \phi \, |^4 
       + |\, \partial_x \phi \, |^2 \, |\, \partial_x U \, |^2
       + |\, \partial_x^2 \phi \, |^2 
       + |\, \phi \, |^2 \,\big( \,|\, \partial_x U  \, |^4
       +  |\, \partial_x^2 U \, |^2
       \,\big) \,\Big)
       \, \mathrm{d}x. 
\end{aligned}
\end{align}
By using Lemma 2.2, Lemma 4.1, Proposition 4.1, and the Sobolev inequality,  
each term in the second term on the right-hand side of (6.9) 
is estimated as follows:
\begin{align}
\begin{aligned}
&\int  
    <\partial_x \phi>^{p-1}
    |\, \partial_x \phi \, |^4 \, \mathrm{d}x
\le
\| \, \partial_x \phi\, \|_{L^{\infty}}^2\int
<\partial_x \phi>^{p-1}
    |\, \partial_x \phi \, |^2 \, \mathrm{d}x
\\[5pt]
&\quad \leq C_0\, 
      \| \, \partial_x \phi \, \|_{L^2} \, 
      \| \, \partial_x^2 \phi \, \|_{L^2}  
      \ll \partial_x\phi \gg _{\infty}^{1-p} 
\\[5pt]
&\quad\leq C_0\, \ll \partial_x\phi \gg _{\infty}^{\frac{5}{2}\, (1-p)}
\\
&\quad \quad\times 
      \Big( \, 
      \int  
       <\partial_x \phi>^{p-1} 
     |\, \partial_x \phi \, |^2 
     \, \mathrm{d}x
      \, \Big)^{\frac{1}{2}} 
      \Big( \, 
      \int 
       <\partial_x \phi>^{2 (p-1)} 
     |\, \partial_x^2 \phi \, |^2 
     \, \mathrm{d}x
      \, \Big)^{\frac{1}{2}};
\end{aligned}
\end{align}
\begin{align}
% \begin{aligned}
\int 
    <\partial_x \phi>^{p-1}
    |\, \partial_x \phi \, |^2 \, |\, \partial_x U \, |^2 \, \mathrm{d}x
\leq C\, 
      \int 
       <\partial_x \phi >^{p-1} 
     |\, \partial_x \phi \, |^2 
     \, \mathrm{d}x\ ;
% \end{aligned}
\end{align}
\begin{align}
% \begin{aligned}
\int 
    <\partial_x \phi>^{p-1}
    |\, \partial_x^2 \phi \, |^2 
    \, \mathrm{d}x
\leq \ \ll \partial_x\phi \gg _{\infty}^{1-p}
      \int  
       <\partial_x \phi>^{2 \, (p-1)} 
     |\, \partial_x^2 \phi \, |^2 
     \, \mathrm{d}x\ ;
% \end{aligned}
\end{align}
\begin{align}
\begin{aligned}
& \int 
    <\partial_x \phi >^{p-1}
    |\, \phi \, |^2 \, (|\, \partial_x U  \, |^4+|\, \partial_x^2 U  \, |^2)
    \, \mathrm{d}x \\[5pt]
& \quad \le
\int 
    |\, \phi \, |^2 \, \big(\, |\, \partial_x U  \, |^4+|\, \partial_x^2 U  \, |^2\, \big)
    \, \mathrm{d}x 
\\[5pt]
& \quad %\quad 
\leq C_0 \, \big( \, \| \, \partial_x U\, \|_{L^{\infty}}^4
+ \| \, \partial_x^2 U \, \|_{L^{\infty}}^2 \, \big)
\leq C_0(1+t)^{-2}.
\end{aligned}
\end{align}
% \begin{equation}
% \int 
%     <\partial_x \phi>^{p-1}
%     |\, \phi \, |^2 \, |\, \partial_x^2 U  \, |^2
%     \, \mathrm{d}x
% \leq C_0 \, \| \, \partial_x^2 U \, \|_{L^{\infty}}^2
% \leq C_0(1+t)^{-2}. 
% \end{equation}
Similarly, each term in $I_{2}$ 
is estimated as follows:
\begin{align}
\begin{aligned}
&\int  
    |\, \partial_x \phi \, |^3 \, 
    |\, \partial _t \partial_x U \, | \, \mathrm{d}x\\
&\leq 2 \, 
      \| \, \partial_x \phi \, \|_{L^2}^2 \, 
      \| \, \partial_x^2 \phi \, \|_{L^2} \, 
    \|\, \partial _t \partial_x U \,\|_{L^2}
\\[5pt]
%&\leq 2 \, \ll \, \partial_x\phi \, \gg _{L^{\infty}_{t, \, x}}^{\frac{3}{2}\, (1-p)}
%      \left( \, 
%      \int ^{\infty }_{-\infty } 
%       <\, \partial_x \phi \, >^{p-1} 
%     |\, \partial_x \phi \, |^2 
%     \, \mathrm{d}x
%       \, \right)^{\frac{1}{2}} \\
%& \quad \times
%      \left( \, 
%      \int ^{\infty }_{-\infty } 
%       <\, \partial_x \phi \, >^{2 \, (p-1)} 
%     |\, \partial_x^2 \phi \, |^2 
%     \, \mathrm{d}x
%      \, \right)^{\frac{1}{2}} \, 
%      \int ^{\infty }_{-\infty } 
%       |\, \partial_x \phi \, | \, 
%    |\, \partial _t \partial_x U \, |
%     \, \mathrm{d}x\\
&\leq C \ll \partial_x\phi \gg _{\infty}^{2\, (1-p)}
\\
&\quad \times
      \Big(\, 
      \int  
       <\partial_x \phi>^{p-1} 
     |\, \partial_x \phi \, |^2 
     \, \mathrm{d}x\, 
      \Big)\, \Big( \, 
      \int  
       <\partial_x \phi>^{2 \, (p-1)} 
     |\, \partial_x^2 \phi \, |^2 
     \, \mathrm{d}x
      \, \Big)^{\frac{1}{2}} \\
&\leq C \, \ll \partial_x\phi \gg _{\infty}^{\frac{5}{2}\, (1-p)}
\\
&\quad \times
      \Big(\, 
      \int  
       <\partial_x \phi>^{p-1} 
     |\, \partial_x \phi \, |^2 
     \, \mathrm{d}x\, 
      \Big)^{\frac 1 2}
\Big( \, 
      \int  
       <\partial_x \phi>^{2 \, (p-1)} 
     |\, \partial_x^2 \phi \, |^2 
     \, \mathrm{d}x
      \, \Big)^{\frac{1}{2}}\ ;
% & \quad \times
%       \int ^{\infty }_{-\infty } 
%       <\, \partial_x \phi \, >^{\frac{1}{2} \, (p-1)} 
%        |\, \partial_x \phi \, | \, 
%     |\, \partial _t \partial_x U \, |
%      \, \mathrm{d}x\\
% &\leq 2 \, \| \, \partial _t \partial_x U(t) \, \|_{L^2} \, 
%       \ll \, \partial_x\phi \, \gg _{L^{\infty}_{t, \, x}}^{2\, (1-p)}
%       \left\| \, \int ^{\infty }_{-\infty } 
%        <\, \partial_x \phi \, >^{p-1} 
%      |\, \partial_x \phi \, |^2 
%      \, \mathrm{d}x \, \right\|_{L^{\infty}_{t}}^{\frac{1}{2}}\\
% & \times 
%       \left( \, 
%       \int ^{\infty }_{-\infty } 
%        <\, \partial_x \phi \, >^{p-1} 
%      |\, \partial_x \phi \, |^2 
%      \, \mathrm{d}x
%        \, \right)^{\frac{1}{2}} 
%       \left( \, 
%       \int ^{\infty }_{-\infty } 
%        <\, \partial_x \phi \, >^{2 \, (p-1)} 
%      |\, \partial_x^2 \phi \, |^2 
%      \, \mathrm{d}x
%       \, \right)^{\frac{1}{2}},
\end{aligned}
\end{align}
\begin{align}
\begin{aligned}
&
\int  
    |\, \partial_x \phi \, |^2 \, 
    |\, \partial _t \partial_x U \, | \, 
    |\,\partial_x U\,| 
    \, \mathrm{d}x\\
%  \leq \| \, \partial _t \partial_x U(t) \, \|_{L^{\infty}} \, 
%       \| \, \partial_x U(t) \, \|_{L^{\infty}} \\
&\qquad 
 \leq C\, \ll \partial_x \phi \gg_{\infty}^{1-p}
% \\
%  & \times     
       \int 
        <\partial_x \phi>^{p-1} 
      |\, \partial_x \phi \, |^2 \, \mathrm{d}x \ ;
\end{aligned}
\end{align}
\begin{align}
\begin{aligned}
&\int 
    |\, \partial_x \phi \, | \, 
    |\, \partial_x^2 \phi \, | \, 
    |\, \partial _t \partial_x U \, | \, 
    |\, \partial_x U\,| 
    \, \mathrm{d}x\\
&\leq C\, 
      \left( \,  
      \| \, \partial_x \phi \, \|_{L^{2}}^2 
      + \| \, \partial_x^2 \phi \, \|_{L^{2}}^2 
      \, \right)\\
&\leq C\, 
%      \| \, \partial _t \partial_x U(t) \, \|_{L^{\infty}} \, 
      %\left( \,  
      \ll \partial_x\phi \gg _{\infty}^{1-p}
      \int 
       <\partial_x \phi>^{p-1} 
     |\, \partial_x \phi \, |^2 
     \, \mathrm{d}x \\% \right. \\
& \quad %\left. 
     + C\, 
%      \| \, \partial _t \partial_x U(t) \, \|_{L^{\infty}} \, 
      \ll \partial_x\phi \gg _{\infty}^{2\,(1-p)}
      \int 
       <\, \partial_x \phi \, >^{2\,(p-1)} 
     |\, \partial_x^2 \phi \, |^2 
     \, \mathrm{d}x \ ; %
%     \, \right),
\end{aligned}
\end{align}
\begin{align}
\begin{aligned}
&\int 
    |\, \phi \, | \, 
    |\, \partial_x \phi \, | \, 
    |\, \partial _t \partial_x U \, | \, 
    \big(\, |\, \partial_x U \, |^2+|\, \partial_x^2 U \, |\, \big)
    \, \mathrm{d}x\\
&\quad\leq C \, 
      \int  
       \Big( \,  
      |\, \partial_x \phi \, |^2
     + |\, \phi \, |^2 \, 
    |\, \partial _t \partial_x U \, |^2 \big( \, 
    \, 
    |\, \partial_x U \, |^4 +|\, \partial_x^2 U \, |^2 \, \big) \, \Big) 
    \, \mathrm{d}x\\
% &\leq \frac{1}{2} \, 
%      \| \, \partial _t \partial_x U(t) \, \|_{L^{\infty}} \, 
%      \| \, \partial_x U(t) \, \|_{L^{\infty}} \, 
%      \int ^{\infty }_{-\infty } 
%       \phi^2 \,  \partial_x U 
%      \, \mathrm{d}x \\% \right. \\
&\quad \leq 
      C\, \ll \partial_x\phi \gg _{\infty}^{1-p}
      \int 
       <\partial_x \phi>^{p-1} 
     |\, \partial_x \phi \, |^2 
     \, \mathrm{d}x + C_0(1+t)^{-2}. %
%     \, \right),
\end{aligned}
\end{align}
% \begin{align}
% \begin{aligned}
% &
% \int ^{\infty }_{-\infty } 
%     |\, \phi \, | \, 
%     |\, \partial_x \phi \, | \, 
%     |\, \partial _t \partial_x U \, | \, 
%     |\, \partial_x^2 U \, |
%     \, \mathrm{d}x%\\
%&
% \leq \frac{1}{2} \, 
%      \| \, \partial _t \partial_x U(t) \, \|_{L^{\infty}}^2 \, 
%      \| \, \partial_x U(t) \, \|_{L^{2}}^2 \\
% & \qquad \qquad \qquad \qquad 
%      + \frac{1}{2} \, 
%       \ll \, \partial_x\phi \, \gg _{L^{\infty}_{t, \, x}}^{1-p}
%       \int ^{\infty }_{-\infty } 
%        <\, \partial_x \phi \, >^{p-1} 
%      |\, \partial_x \phi \, |^2 
%      \, \mathrm{d}x, %
%     \, \right),
% \end{aligned}
% \end{align}
% where 
% $ 
% \displaystyle{
% \Vert \, w \, \Vert_{_{L^{\infty}_{t}}}
% :=  \esssup\displaylimits_{t\geq0} \, |\, w(t) \,|}. 
% $
% Noting (6.7) and using the Young inequality, 
% we finally estimate the right-hand side of (6.2) as 
We finally estimate the right-hand side of (6.2) 
by (6.7) and the Young inequality as 
\begin{align}
\begin{aligned}
&\left| \,  
\int 
  \partial _t \big( \, 
  \sigma( \, 
  \partial_x U + \partial_x \phi \,) - 
    \sigma( \, 
    \partial_x U  \, 
    )  \, 
    \big) \, 
  \partial_x^2 \big( \, \sigma( \, \partial_x U  \,)  \, \big)
\, \mathrm{d}x
\, \right|\\
&\quad \leq C \,  \int 
       <\partial_x \phi>^{p-1} 
     |\, \partial _t \partial_x \phi \, |\, 
     \left|\, 
     \partial_x^2 \big( \, \sigma ( \, \partial_x U  \,)  \, \big) 
     \, \right|
     \, \mathrm{d}x \\
& \quad \quad
     + C \,  \int 
       |\, \partial_x \phi \, |\, 
       |\, \partial _t \partial_x U \, |\, 
       \left|\, 
     \partial_x^2 \left( \, \sigma\bigl( \, \partial_x U  \, \bigr)  \, \right) 
     \, \right|
     \, \mathrm{d}x. 
\end{aligned}
\end{align}
Each term on the right-hand side of (6.18) is estimated as 
\begin{align}
\begin{aligned}
& \int 
       <\partial_x \phi>^{p-1} 
     |\, \partial _t \partial_x \phi \, |\, 
     \left|\, 
     \partial_x^2 \left( \, \sigma\bigl( \, \partial_x U  \, \bigr)  \, \right) 
     \, \right|
     \, \mathrm{d}x \\
&\quad \leq \epsilon \, \int
       <\, \partial_x \phi \, >^{p-1} 
     |\, \partial _t \partial_x \phi \, |^2 
     \, \mathrm{d}x
%& 
+ C_{\epsilon} \,  \left\|\, 
     \partial_x^2 \left( \, \sigma\bigl( \, \partial_x U  \, \bigr)  \, \right) (t)
     \, \right\|_{L^2}^2
\\
&\quad \leq \epsilon \, \int
       <\, \partial_x \phi \, >^{p-1} 
     |\, \partial _t \partial_x \phi \, |^2 
     \, \mathrm{d}x
+ C_{\epsilon}\,(1+t)^{-2}, 
\end{aligned}
\end{align}
and
\begin{align}
\begin{aligned}
& \int 
       |\, \partial_x \phi \, |\, 
       |\, \partial _t \partial_x U \, |\, 
       \left|\, 
     \partial_x^2 \big( \, \sigma( \, \partial_x U  \,)  \, \big) 
     \, \right|
     \, \mathrm{d}x\\
% &\quad \leq C \,  
%      \| \, \partial _t \partial_x U(t) \, \|_{L^{\infty}}^2 \, 
%      \left\|\, 
%      \partial_x^2 \left( \, \sigma\bigl( \, \partial_x U  \, \bigr)  \, \right) (t)
%      \, \right\|_{L^2}^2\\
& \quad 
\le C \, \ll \partial_x\phi \gg _{\infty}^{1-p}
      \int 
       <\partial_x \phi>^{p-1} 
     |\, \partial_x \phi \, |^2 
     \, \mathrm{d}x + C\,(1+t)^{-2}. 
\end{aligned}
\end{align}
Then, substituting all the estimates (6.5)$\sim$(6.20) into (6.2),
choosing $\epsilon$ suitably small, and integrating the 
resultant formula with respect to $t$, we arrive at
\begin{align}
\begin{aligned}
&
\int 
\big| \, \partial _x \big( \, 
  \sigma( \, 
  \partial_x U + \partial_x \phi \,) - 
    \sigma( \, 
    \partial_x U  \,)  \, 
    \big)
    \, \big|^2 \, \mathrm{d}x
% & + \int ^{\infty }_{-\infty } 
%   \Big( \, 
%   \sigma'\bigl( \, 
%   \partial_x U + \partial_x \phi \, \bigr) - 
%     \sigma'\bigl( \, 
%     \partial_x U  \, 
%     \bigr)  \, 
%     \Big)^2 \, 
%     ( \, \partial_x^2 U \, )^2  
%     \, \mathrm{d}x\\
\\ &%\quad 
+
\int^{t}_{0} \int 
     <\partial_x \phi>^{p-1} 
|\, \partial_t \partial_x \phi \, |^2 
\, \mathrm{d}x\mathrm{d}\tau 
\leq C_0\ll \partial_x\phi \gg _{\infty}^{3\,(1-p)}
\quad \big(\,t \in [\,0,\:T\,] \, \big),
% \\
% & \quad
% \leq C_0
% \, \left( \, 1 
% + \ll \, \partial_x\phi \, \gg _{L^{\infty}_{t, \, x}}^{1-p}
% + \ll \, \partial_x\phi \, \gg _{L^{\infty}_{t, \, x}}^{2\,(1-p)} 
% \ll \partial_x\phi \gg _{\infty}^{3(1-p)}\\ 
% \, \right)\\
% &\quad \quad
% + C_0 \ll \partial_x\phi \gg _{\infty}^{\frac{5}{2}\, (1-p)} 
%      \int^{t}_{0} 
%       \Big( \, 
%       \int 
%        <\partial_x \phi>^{p-1} 
%      |\, \partial_x \phi \, |^2 
%      \, \mathrm{d}x
%        \, \Big)^{\frac{1}{2}} \\
% &\quad\quad 
%       \times 
%       \left( \, 
%       \int 
%        <\, \partial_x \phi \, >^{2 \, (p-1)} 
%      |\, \partial_x^2 \phi \, |^2 
%      \, \mathrm{d}x
%       \, \right)^{\frac{1}{2}} 
%       \, \mathrm{d}\tau. 
\end{aligned}
\end{align}
where we used the estimate
\begin{align}
\begin{aligned}
&\int^{t}_{0} 
      \Big(\,
      \int 
       <\partial_x \phi>^{p-1} 
     |\, \partial_x \phi \, |^2 
     \, \mathrm{d}x\mathrm{d}\tau \,
       \Big)^{\frac{1}{2}}\,
      \Big(\,
      \int  
       <\partial_x \phi>^{2 (p-1)} 
     |\, \partial_x^2 \phi \, |^2 
     \, \mathrm{d}x \, 
      \,\Big)^{\frac{1}{2}} 
      \, \mathrm{d}\tau \\
&\quad \leq
\Big( \,\int^{t}_{0} 
      \int 
       <\partial_x \phi>^{p-1} 
     |\, \partial_x \phi \, |^2 
     \, \mathrm{d}x\mathrm{d}\tau
      \,\Big)^{\frac{1}{2}}\\
&\quad\quad \times
      \big( \,\int^{t}_{0} 
      \int 
       <\partial_x \phi>^{2 (p-1)} 
     |\, \partial_x^2 \phi \, |^2 
     \, \mathrm{d}x\mathrm{d}\tau
      \, \Big)^{\frac{1}{2}}
\leq C_0  \,\ll \partial_x\phi \gg _{\infty}^{\frac{1}{2} \, (1-p)}. 
\end{aligned}
\end{align}
% Therefore, substituting (6.29) into (6.28), 
% we have 
% \begin{align}
% \begin{aligned}
% & \int ^{\infty }_{-\infty } 
%      <\, \partial_x \phi \, >^{2\,(p-1)} 
% |\, \partial_x^2 \phi \, |^2 
% \, \mathrm{d}x\\
% & + \int ^{\infty }_{-\infty } 
%   \Big( \, 
%   \sigma'\bigl( \, 
%   \partial_x U + \partial_x \phi \, \bigr) - 
%     \sigma'\bigl( \, 
%     \partial_x U  \, 
%     \bigr)  \, 
%     \Big)^2 \, 
%     ( \, \partial_x^2 U \, )^2  
%     \, \mathrm{d}x\\
%    &+ \int^{t}_{0} \int ^{\infty }_{-\infty } 
%      <\, \partial_x \phi \, >^{p-1} 
% |\, \partial_t \partial_x \phi \, |^2 
% \, \mathrm{d}x\mathrm{d}\tau \\
% &\leq C ( \, \phi_{0}, \: \partial_x\phi_{0}, \: \partial_x^2\phi_{0} \, )
% \, \left( \, 1 
% + \ll \, \partial_x\phi \, \gg _{L^{\infty}_{t, \, x}}^{3\,(1-p)} 
% \, \right). 
% \end{aligned}
% \end{align}
Finally, if we note the estimates (5.7) and (5.8) imply
$$
\int 
     <\partial_x \phi>^{2\,(p-1)} 
     |\, \partial_x^2 \phi \, |^2 
     \, \mathrm{d}x
\le  
C\int  
  \big|\, \partial _x \big( \, 
  \sigma ( \, 
  \partial_x U + \partial_x \phi \,) - 
    \sigma ( \, 
    \partial_x U  \, 
    )  \, 
    \big)\big|^2  
    \, \mathrm{d}x+C,
$$
the estimate (6.21) immediately implies the desired 
{\it a priori} estimate for $\partial_x^2\phi$.
Thus the proof of Proposition 6.1 is completed.

\medskip

Now, combining Proposition 4.1, Proposition 5.1, and Proposition 6.1,
we show the the following uniform boundedness of $\|\partial_x\phi\|_{L^{\infty}}$
which plays the essential role to control the nonlinearity of $\sigma$.
The proof is motivated by an idea in Kanel' [16].

% To complete the proof of Proposition 6.1, 
% we claim the following uniform boundedness of $\partial_x\phi$. 

\medskip

\noindent
{\bf Lemma 6.1.} %(uniform boundedness for $\partial_x \phi$){\bf .}
\quad {\it For $3/7 <p <1$,
there exists a positive constant $C_0$ such that 
$$
% \displaystyle{ 
% \sup_{t\geq0, \; x \in \mathbb{R}} \, |\, \partial_x \phi(t, \,x) \,| \leq C.
% }
\|\,\partial_x\phi (t)\,\|_{L^{\infty}} \le C_0 \qquad \big(\,t \in [\,0,\:T\,] \, \big).
$$
}

% \medskip

% \noindent
% In fact, 
% once Lemma 6.1 holds true, we immediately have Proposition 6.1 from (6.28). 
% Therefore, we show Lemma 6.1 
% with the aid of the Sobolev type inequality 
% motivated by an idea in Kanel' (\cite{kanel}). 

\medskip

{\bf Proof of Lemma 6.1}.
By the Schwarz inequality, we have 
for $a>0$
\begin{align}
\begin{aligned}
&<\partial_x \phi(t,\, x)>^a 
= 1+ \int_{-\infty}^x\frac{\partial}{\partial y}<\partial_x\phi(t,y)>^a\,dy
\\
&\leq 1+ a
\int <\partial_x\phi>^{a-2}
     |\, \partial_x\phi\, | \, 
     |\, \partial_x^2\phi \, |
     \, \mathrm{d}x\\
&
\leq 1+ a\, \Big(\, 
      \int 
       <\partial_x\phi>^{p-1} 
     |\, \partial_x\phi \, |^2 
     \, \mathrm{d}x
      \, \Big)^{\frac{1}{2}} \,
      \Big( \, 
      \int 
       <\partial_x\phi>^{2a-3-p} 
     |\, \partial_x^2 \phi \, |^2 
     \, \mathrm{d}x
     \,  \Big)^{\frac{1}{2}}. 
\end{aligned}
\end{align}
If we choose $a=(3\, p+1)/2$, (6.23) gives
\begin{align}
\begin{aligned}
\ll  \partial_x \phi \gg_{\infty}^{\frac{3p+1}{2}}
\leq 1 &+ C\, \Big(\,
      \int 
       <\partial_x\phi>^{p-1} 
     |\, \partial_x\phi \, |^2 
     \, \mathrm{d}x
     \, \Big)^{\frac{1}{2}} 
\\
& %\qquad \qquad \qquad \qquad
      \times \big( \,
      \int 
       <\partial_x\phi>^{2\,(p-1)} 
     |\, \partial_x^2 \phi \, |^2 
     \, \mathrm{d}x
      \, \Big)^{\frac{1}{2}}. 
\end{aligned}
\end{align}
which deduces from Proposition 5.1 and Proposition 6.1 that
\begin{equation}
\ll  \partial_x \phi \gg_{\infty}^{\frac{3p+1}{2}}\, \le \,
C_0 \, \ll  \partial_x \phi \gg_{\infty}^{2\,(p-1)}.
\end{equation}
Hence, if we assume
$$
\frac{3\, p+1}{2} < 2\, (p-1)\qquad \bigg(\,\Leftrightarrow p > \frac{3}{7}\,\bigg),
$$
we obtain for $3/7<p<1$
\begin{equation}
\ll \partial_x\phi \gg _{\infty}
\leq \,C_0. 
\end{equation}
Thus, the proof of Lemma 6.1 is completed. 

\medskip

By Proposition 4.1, Proposition 5.1, and Proposition 6.1 with the aid of
Lemma 6.1, we obtain the energy estimate
\begin{align}
\begin{aligned}
\| \, \phi(t) \, \|_{{H}^{2}}^{2} 
& + \int^{t}_{0} 
   \big\| \, 
   (\, \sqrt{\, \partial_x U} \:  \phi \,)(\tau) 
   \, \big\|_{L^2}^{2} \, \mathrm{d}\tau 
\\
& %\quad 
+ \int^{t}_{0}
   \big( 
 \left\| \, \partial_{x}\phi(\tau) \, \right\|_{{H}^1}^{2} 
 + \left\| \, \partial_{t}\partial_{x}\phi(\tau) \, \right\|_{L^2}^{2} 
   \, \big) \, \mathrm{d}\tau 
\leq C_0
\qquad \big( \, t \in [\, 0, \: T\, ] \,\big).
\end{aligned}
\end{align}
Therefore, in order to accomplish the proof of Theorem 3.3, 
it suffices to show the following {\it a priori} estimate: 
\begin{equation}
\int_0^t \| \, \partial_{x}^3\phi(\tau) \,\|_{L^{2}}^2\, \mathrm{d}\tau 
\le C_0 \qquad \big( \, t \in [\, 0, \: T\, ] \,\big).
\end{equation}
The estimate (6.28) is directly obtained by
the equation (6.1) and the estimate (6.27) as follows. 
The equation (6.1) is rewritten as
\begin{align*}
\begin{aligned}
     \sigma'( \, \partial_x U &+ \partial_x \phi \, ) 
     \partial_x^3\phi = \partial_t\partial_x\phi+ \partial_x^2 \big( \, f(U+\phi) - f(U) \, \big) 
\\[5pt]
  &-\sigma''( \, \partial_x U + \partial_x \phi \, )|\, \partial_x U + \partial_x \phi \,|^2 
    - \sigma' (\, \partial_x U + \partial_x \phi \, )\partial_x^3 U. 
\end{aligned}
\end{align*}
Then, by the estimate (6.27), Lemma 2.2, and the Sobolev inequality, 
it holds
\begin{align*}
\begin{aligned}
\int_0^t \| \, \partial_{x}^3\phi\,\|_{L^{2}}^2\, \mathrm{d}\tau 
&\le C_0 + C_0\, \int_0^t\int | \, \partial_{x}^2\phi\,|^4
\, \mathrm{d}x\mathrm{d}\tau\\[5pt]
&\le C_0 + C_0\, \int_0^t \| \, \partial_{x}^2\phi\,\|_{L^2}^3\,
 \| \, \partial_{x}^3\phi\,\|_{L^2}
\,\mathrm{d}\tau\\[5pt]
&\le C_0 + \frac{1}{2}\, \int_0^t \| \,\partial_{x}^3\phi\,\|_{L^2}^2
\,\mathrm{d}\tau
+ C_0\, \int_0^t \| \,\partial_{x}^2\phi\,\|_{L^2}^2
\,\mathrm{d}\tau,
\end{aligned}
\end{align*}
which implies (6.28).
Thus, the proof of Theorem 3.3 is completed.

\end{document}